\documentclass{amsart}
\usepackage{graphicx}

\newtheorem{theorem}{Theorem}[section]
\newtheorem{lemma}[theorem]{Lemma}

\theoremstyle{definition}

\theoremstyle{remark}
\newtheorem{remark}[theorem]{Remark}

\numberwithin{equation}{section}

\begin{document}

\title[Closed incompressible surfaces of genus two]{Closed incompressible surfaces of genus two in 3-bridge knot complements}

\author{Makoto Ozawa}
\address{Department of Natural Sciences, Faculty of Arts and Sciences, Komazawa University, 1-23-1 Komazawa, Setagaya-ku, Tokyo, 154-8525, Japan}
\email{w3c@komazawa-u.ac.jp}





\begin{abstract}
In this paper, we characterize closed incompressible surfaces of genus two in the complements of 3-bridge knots and links. This characterization includes that of essential 2-string tangle decompositions for 3-bridge knots and links.
\end{abstract}

\maketitle

\section{Introduction and Result}
Possible positions of a knot in the 3-sphere restrict on forms and positions of closed incompressible surfaces in the knot complement, and vice versa.

Allen Hatcher and William Thurston showed that there is no closed incompressible non-peripheral surface in the complements of 2-bridge knots, namely, 2-bridge knots are small (\cite{HT}).
This implies that knots whose complements contain closed incompressible surfaces except for a peripheral torus are not 2-bridge.
Ulrich Oertel classified all closed incompressible non-peripheral surfaces in the complements of Montesinos knots and links by showing that they are obtained from Conway spheres by meridional tubing (\cite{O}).
This implies that if a knot complement contains a closed incompressible surface which is not obtained from Conway spheres by meridional tubing, then the knot is not Montesinos.
William Menasco proved that any closed incompressible surface in the complements of alternating knots and links is meridionally compressible (\cite{M}).
This implies that knots whose complements contain closed incompressible and meridionally incompressible surfaces are not alternating.
The author showed that any closed incompressible surfaces in the complements of positive knots have non-zero order, that is, there exists a loop in the surface that has non-zero linking number with the knot (\cite{Oz}).
This implies that if the complement of a knot contains a closed incompressible surface with zero order, namely any loop on the surface has zero linking number with the knot, then the knot is not positive.

In this paper, we characterize closed incompressible surfaces of genus two in the complements of 3-bridge knots and links. This characterization includes that of essential 2-string tangle decompositions for 3-bridge knots and links.
As above observation, we can determine whether a given knot is 3-bridge when the knot complement contains a closed incompressible surface of genus two.

Hereafter, let $K$ be a knot or link in the 3-sphere $S^3$, and $F$ a closed surface embedded in $S^3$ which is disjoint from $K$ or intersects $K$ transversely.

We say that a disk $D$ embedded in $S^3-K$ is a {\em compressing disk} for $F-K$ if $D\cap F=\partial D$ and $\partial D$ is essential in $F-K$.
If there exists a compressing disk $D$ for $F-K$, by cutting $F$ along $\partial D$ and pasting two parallel copies of $D$, we obtain a new closed surface $F'$.
Such an operation is called a {\em compression}.
We say that $F$ is {\em compressible} in $S^3-K$ if there exists a compressing disk for $F-K$, and that $F$ is {\em incompressible} if $F$ is not compressible.

Next, 
we say that a disk $D$ embedded in $S^3$ is a {\em meridional compressing disk} for $F$ if $D$ intersects $K$ in one point in the interior of $D$, $D\cap F=\partial D$ and $\partial D$ is not parallel to a component of $\partial N(K;F)$ in $F$.
If there exists a meridional compressing disk $D$ for $F$, by cutting $F$ along $\partial D$ and pasting two parallel copies of $D$, we obtain a new closed surface $F'$.
Such an operation is called a {\em meridional compression}.
We say that $F$ is {\em meridionally compressible} in $S^3-K$ if there exists a meridional compressing disk for $F$, and that $F$ is {\em meridionally incompressible} if $F$ is not meridionally compressible.

It is fundamental and important that any closed surface becomes to be incompressible by some compressions, and any closed incompressible surface becomes to be meridionally incompressible by some meridional compressions.
Conversely, any closed incompressible surface can be obtained from closed incompressible and meridionally incompressible surfaces by some meridional tubings, and any closed surface can be obtained from closed incompressible surfaces by some tubings.
Hence, closed incompressible and meridionally incompressible surfaces are basic among all closed surfaces.

In the case that $F$ is a genus two closed incompressible surface in $S^3-K$, there are three types of surfaces obtained by meridional compressions.

\begin{description}
	\item[Type I] a genus two closed surface which is disjoint from $K$
	\item[Type II] a genus one closed surface which intersects $K$ in two points
	\item[Type III] a genus zero closed surface which intersects $K$ in four points, that is, an essential Conway sphere
\end{description}

\begin{figure}[htbp]
	\begin{center}
	\begin{tabular}{ccc}
		\includegraphics[trim=0mm 0mm 0mm 0mm, width=.2\linewidth]{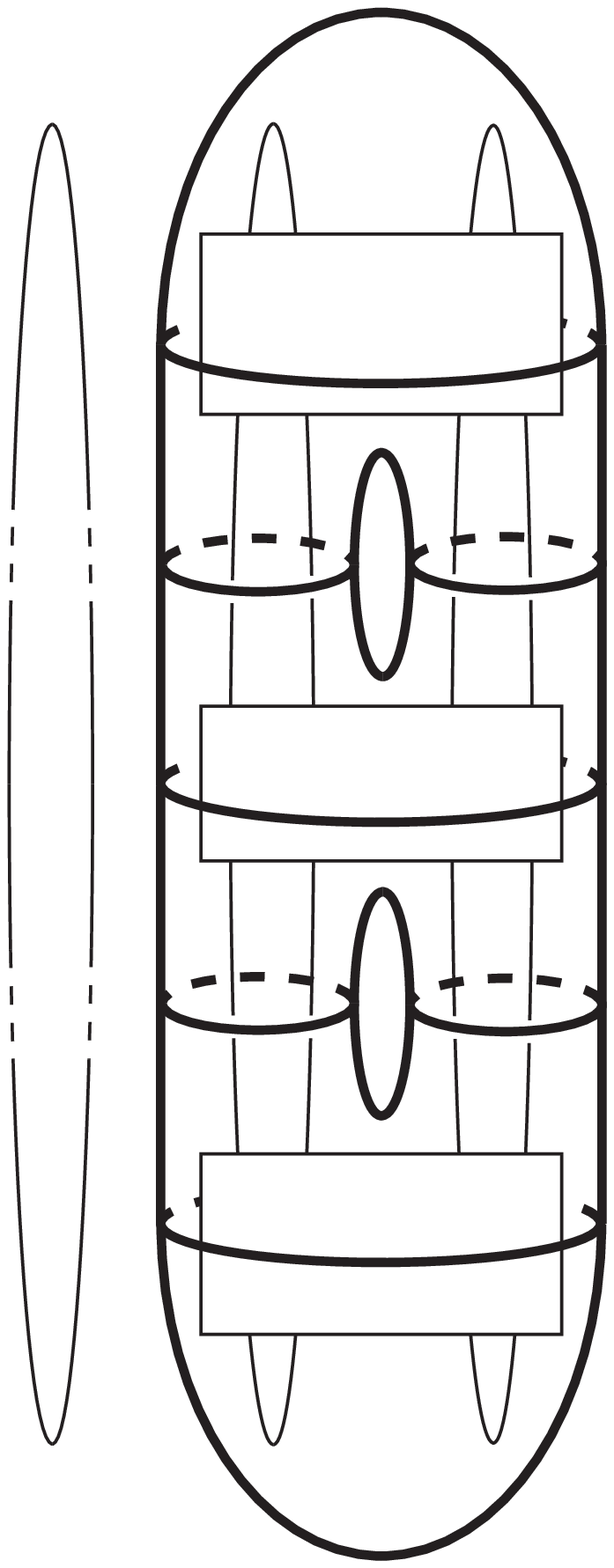} &
		\includegraphics[trim=0mm 0mm 0mm 0mm, width=.25\linewidth]{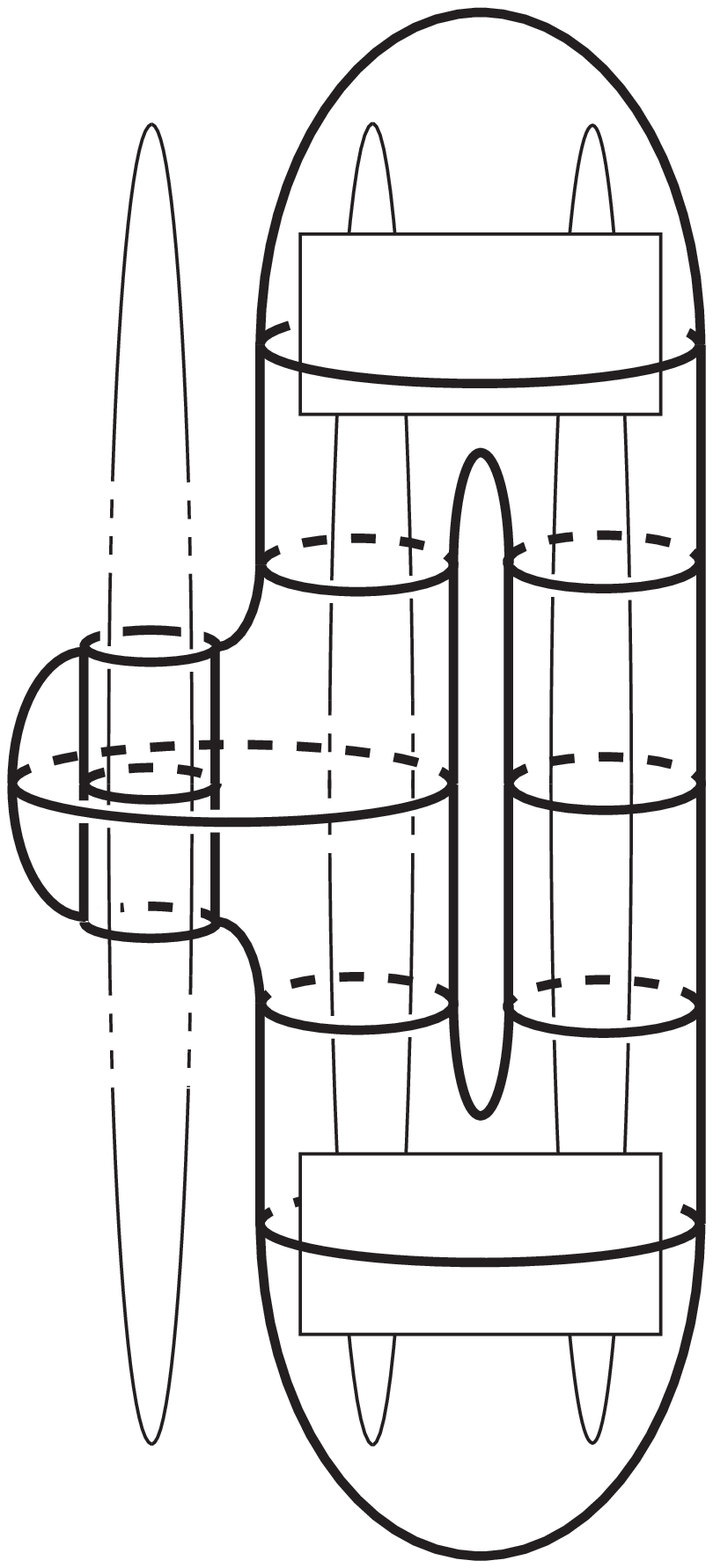} &
		\includegraphics[trim=0mm 0mm 0mm 0mm, width=.3\linewidth]{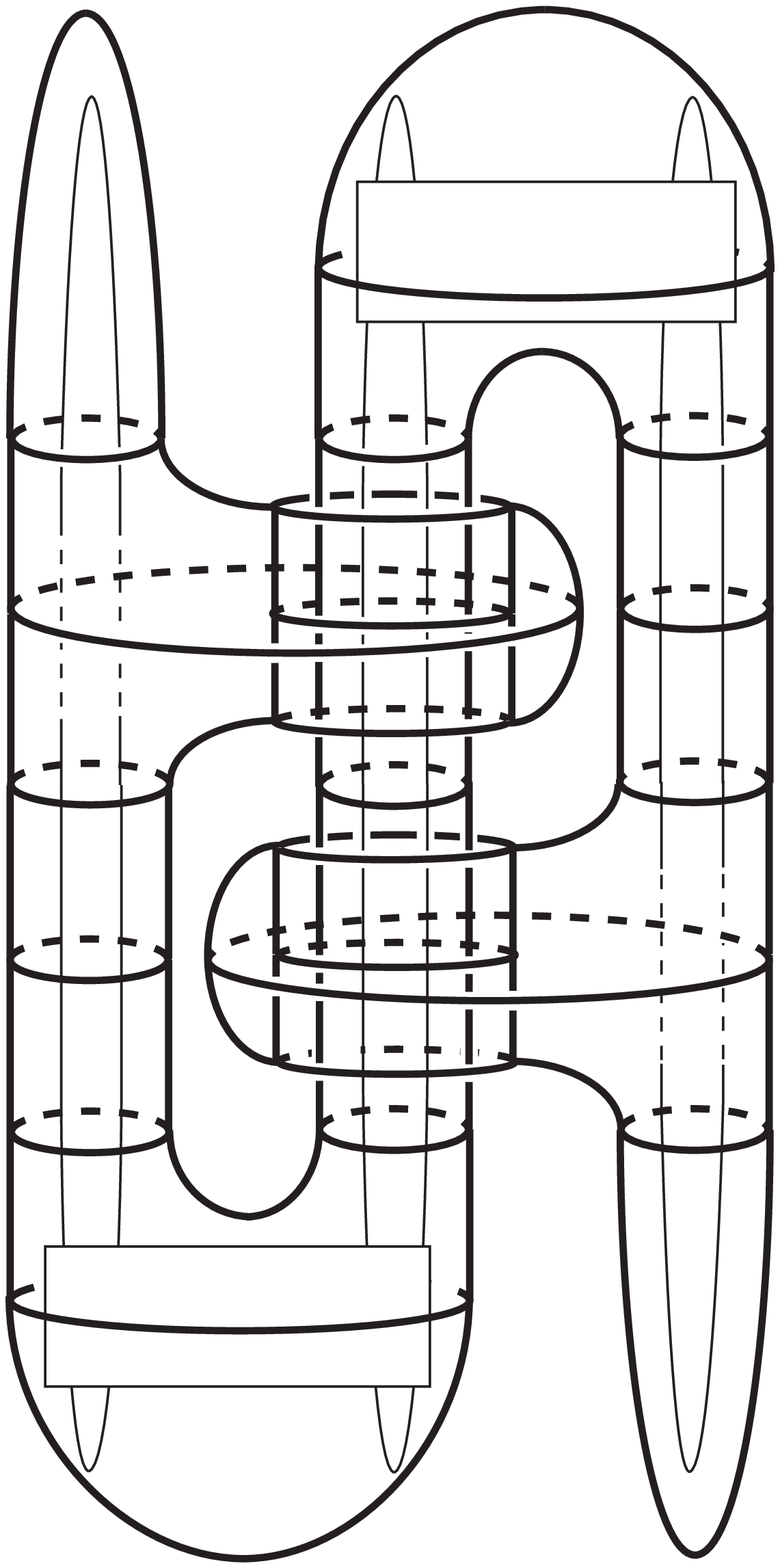} \\
		I-a & I-b & I-c
	\end{tabular}
	\end{center}
	\caption{genus two closed surfaces of Type I}
	\label{figure I}
\end{figure}

\begin{figure}[htbp]
	\begin{center}
	\begin{tabular}{ccc}
		\includegraphics[trim=0mm 0mm 0mm 0mm, width=.25\linewidth]{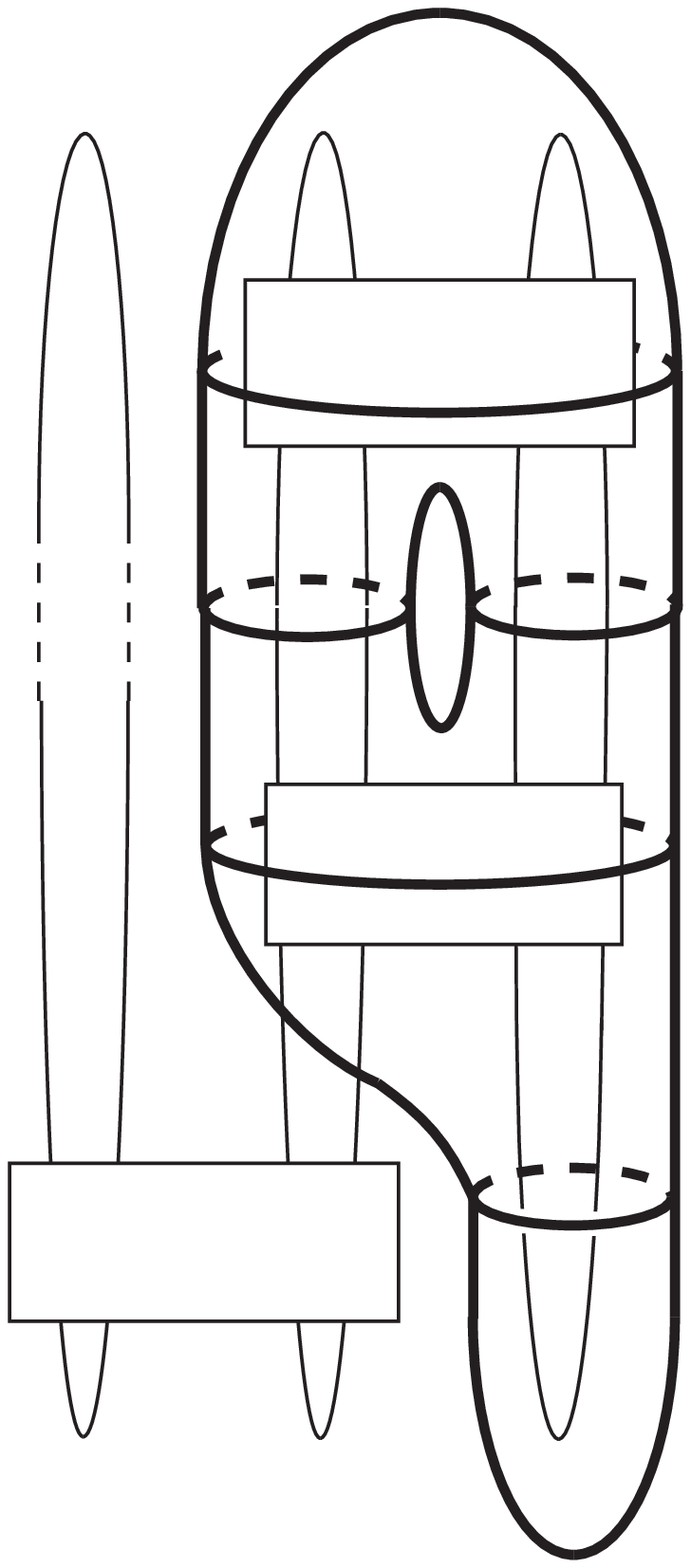} &
		\includegraphics[trim=0mm 0mm 0mm 0mm, width=.23\linewidth]{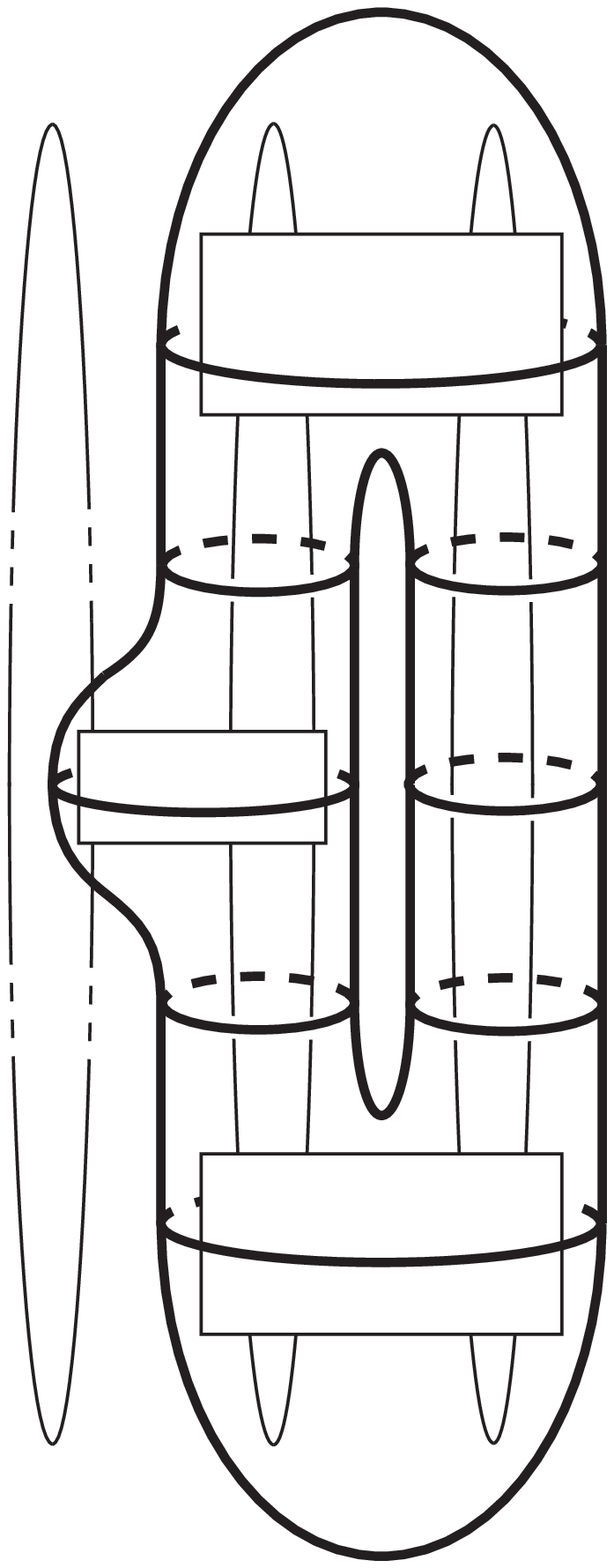} &
		\includegraphics[trim=0mm 0mm 0mm 0mm, width=.27\linewidth]{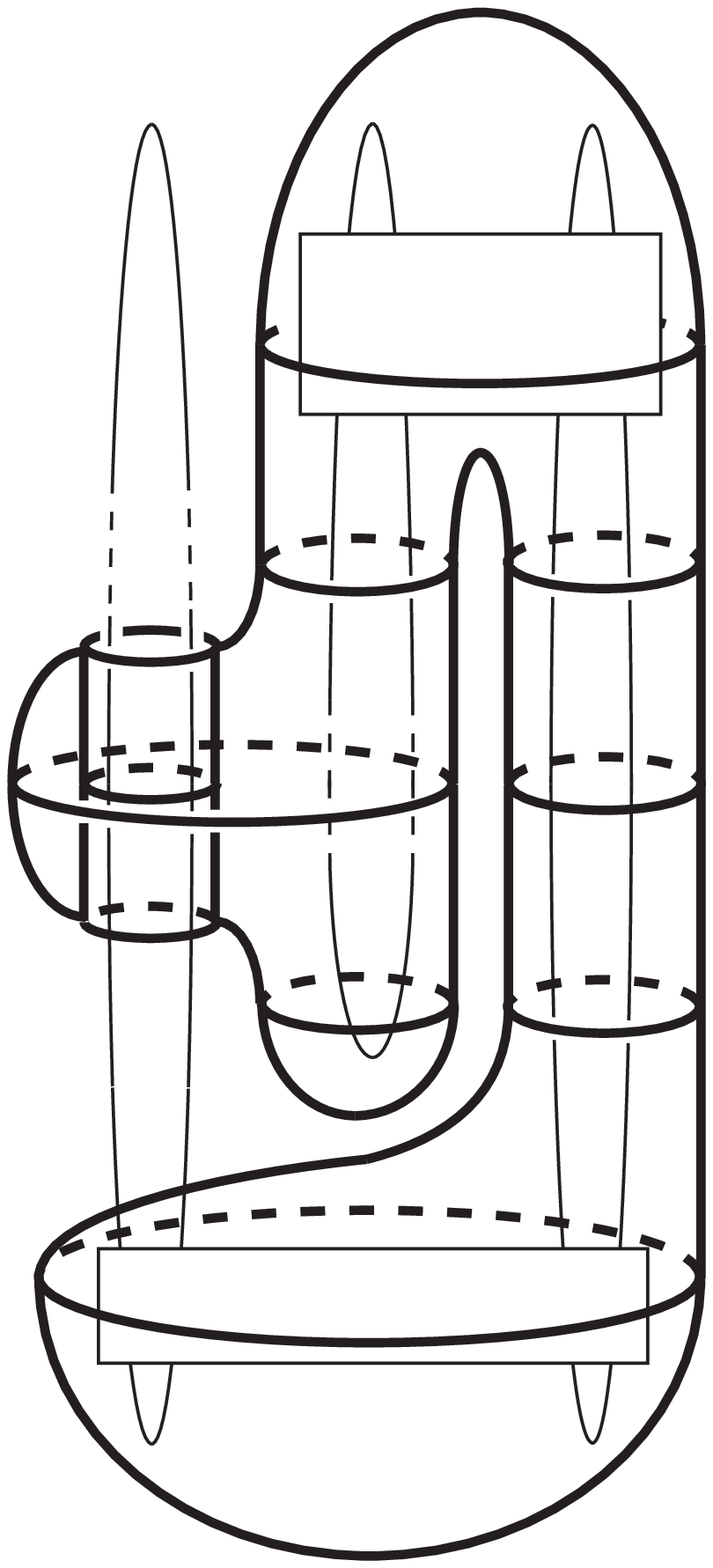} \\
		II-a & II-b & II-c
	\end{tabular}
	\end{center}
	\caption{genus one closed surfaces of Type II}
	\label{figure II}
\end{figure}

\begin{figure}[htbp]
	\begin{center}
	\begin{tabular}{cc}
		\includegraphics[trim=0mm 0mm 0mm 0mm, width=.25\linewidth]{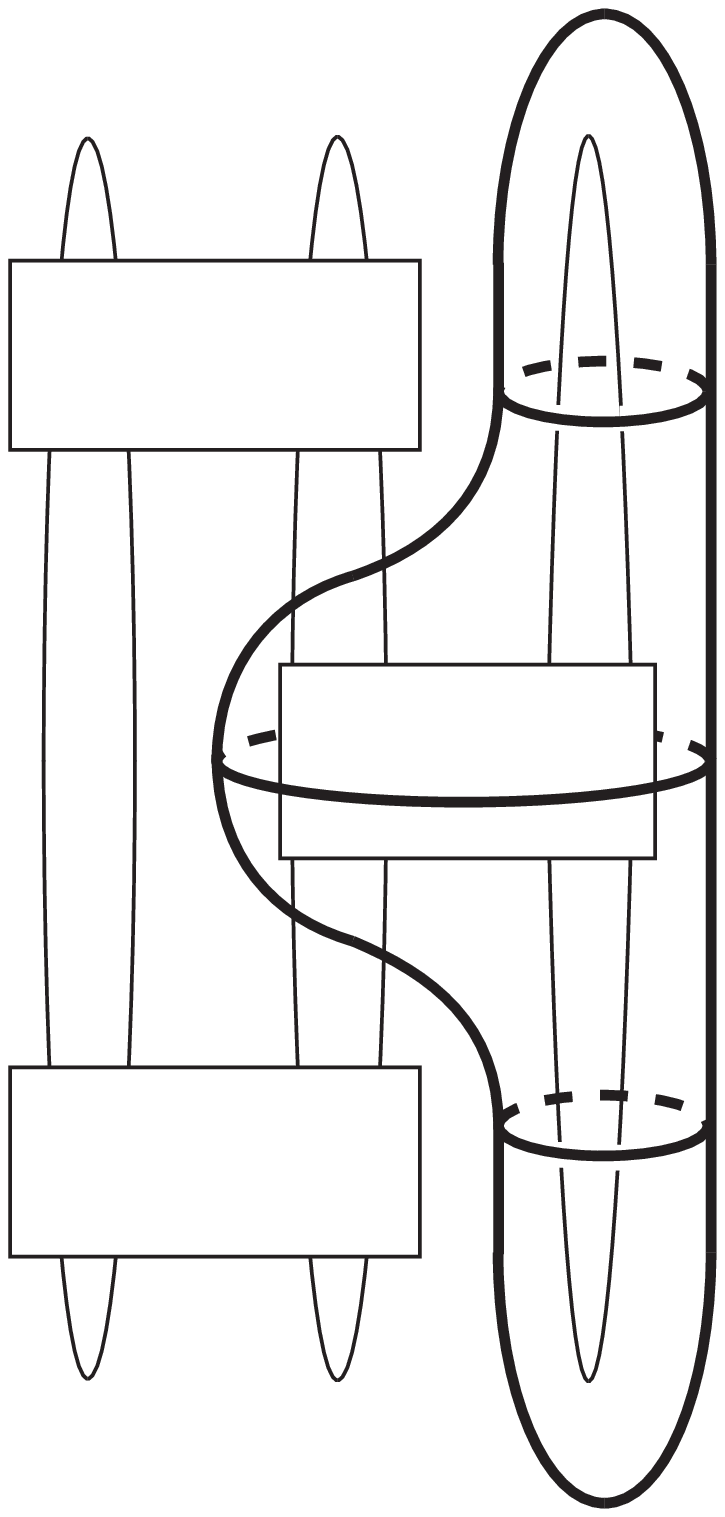} &
		\includegraphics[trim=0mm 0mm 0mm 0mm, width=.25\linewidth]{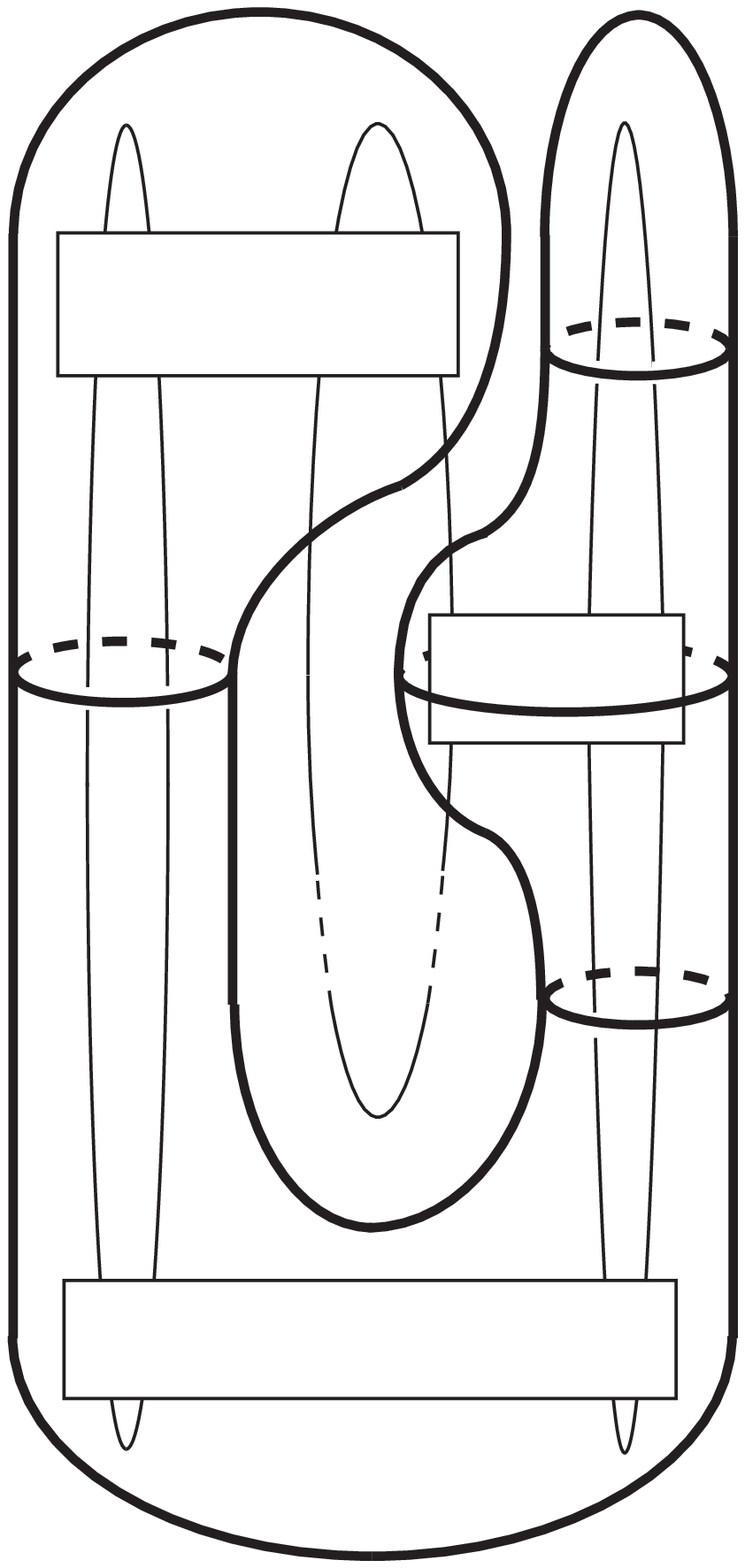} \\
		III-a & III-b
	\end{tabular}
	\end{center}
	\caption{genus zero closed surfaces of Type III}
	\label{figure III}
\end{figure}

In the following three theorems, we characterize closed incompressible surfaces of genus two in the 3-bridge knot or link complement.
Figures \ref{figure I}, \ref{figure II} and \ref{figure III} illustrate closed incompressible and meridionally incompressible surfaces of Type I, II and III respectively.
Here, all boxes indicate braids and dotted lines indicate monotone arcs.
A description of ``types for triples'' and a definition of a ``gradient disk'' are given in the next section.

\begin{theorem}\label{I}
Let $K$ be a 3-bridge knot or link, and $F$ be a closed incompressible and meridionally incompressible surface of genus two.
Then, the triple $(S^3,K,F)$ is either
\begin{description}
	\item [I-a] a union of two triples of Type $T_0$, or
	\item [I-b] a union of two triples of Type $P_0$, or
	\item [I-c] a union of two triples of Type $A_0\cup T_0$.
\end{description}
Conversely, $F$ is incompressible and meridionally incompressible if there is no gradient disk for all critical points.
\end{theorem}

\begin{theorem}\label{II}
Let $K$ be a 3-bridge knot or link, and $F$ be a closed incompressible and meridionally incompressible surface of genus one which intersects $K$ in two points.
Then, the triple $(S^3,K,F)$ is either
\begin{description}
	\item [II-a] a union of a triple of Type $T_0$ and of Type $D_2$, or
	\item [II-b] a union of two triples of Type $A_1$, or
	\item [II-c] a union of a triple of Type $P_0$ and of Type $D_2\cup A_0$.
\end{description}
Conversely, $F$ is incompressible and meridionally incompressible if there is no gradient disk for all critical points and two points of $K\cap F$.
\end{theorem}

\begin{theorem}\label{III}
Let $K$ be a 3-bridge knot or link, and $F$ be a closed incompressible and meridionally incompressible surface of genus zero which intersects $K$ in four points.
Then, the triple $(S^3,K,F)$ is either
\begin{description}
	\item [III-a] a union of two triples of Type $D_2$, or
	\item [III-b] a union of a triple of Type $D_1\cup D_2$ and of Type $A_1$.
\end{description}
Conversely, $F$ is incompressible and meridionally incompressible if there is no gradient disk for all critical points and four points of $K\cap F$.
\end{theorem}

\section{Preliminary}

\begin{lemma}\label{level}
Let $K$ be a knot or link in a bridge position with respect to the standard Morse function $h:S^3\to \Bbb{R}$, and $F$ be a closed incompressible and meridionally incompressible surface.
Then one of the following holds.
\begin{enumerate}
	\item $K$ is a split link.
	\item $K$ is not thin position.
	\item After an isotopy of $K$ and $F$, there exists a level sphere $S=h^{-1}(t)$ such that each component of $S\cap F$ is essential in both $S-K$ and $F-K$.
\end{enumerate}
\end{lemma}

\begin{proof}
First, we isotope $K$ and $F$ so that $F$ has no inessential saddle, that is, at least one of two level loops appearing in a saddle point bounds a disk in $F$.
See Figure \ref{isotopy}.

\begin{figure}[htbp]
	\begin{center}
		\includegraphics[trim=0mm 0mm 0mm 0mm, width=.6\linewidth]{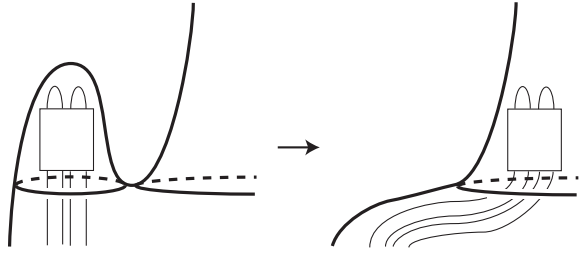}
	\end{center}
	\caption{an isotopy of $K$ and $F$ reducing an inessential saddle}
	\label{isotopy}
\end{figure}

Next, consider the set of level spheres whose intersection with $F$ contains inessential loops.
If the set ranges in all level, then the condition 1 or 2 holds.
Otherwise, the condition 3 holds.
\end{proof}

If a 3-bridge knot or link $K$ is split or composite, then there does not exist a genus two closed incompressible surface.
Hence, for any Type I, II or III, by Lemma \ref{level}, there exists a level sphere $S$ which intersects $F-K$ essentially and decomposes the pair $(S^3,K)$ into two 3-string trivial tangles.
Threrefore, next we concentrate on incompressible and meridionally incompressible surfaces in a 3-string trivial tangle.

\begin{lemma}\label{trivial}
Let $(B,T)$ be an $n$-string trivial tangle and $P$ an incompressible and not $\partial$-parallel surface in $(B,T)$.
Then, one of the following holds.
\begin{enumerate}
	\item $P$ is a disk which is disjoint from $T$ and separates $(B,T)$ into two trivial tangles.
	\item $P$ is a disk which intersects $T$ in one point and separates $(B,T)$ into two trivial tangles.
	\item $P$ is $\partial$-compressible.
\end{enumerate}
\end{lemma}

\begin{proof}
This Lemma can be proved by observing the intersection of $P$ and a system of ``trivializing disks'' for the strings $T$, that is, a union of disjoint $n$ disks $\Delta_1,\ldots,\Delta_n$ such that $\Delta_i\cap T=\partial \Delta_i\cap T=t_i$ and $\Delta_i\cap \partial B=\partial \Delta_i-t_i$, where $t_i$ is a string of $T$.
By cutting and pasting a system of trivializing disks, we assume that the number of the intersection is minimal.
If $P$ does not intersect a system of trivializing disks, then we have the conclusion 1.
If $P$ intersects a system of trivializing disks in an arc, then we have the concludion 2.
Otherwise, we have the conclusion 3.
\end{proof}

Using Lemma \ref{trivial} inductively, we can classify incompressible and meridionally incompressible surfaces in the trivial 3-string tangle as follows.

\begin{figure}[htbp]
	\begin{center}
	\begin{tabular}{ccc}
		\includegraphics[trim=0mm 0mm 0mm 0mm, width=.3\linewidth]{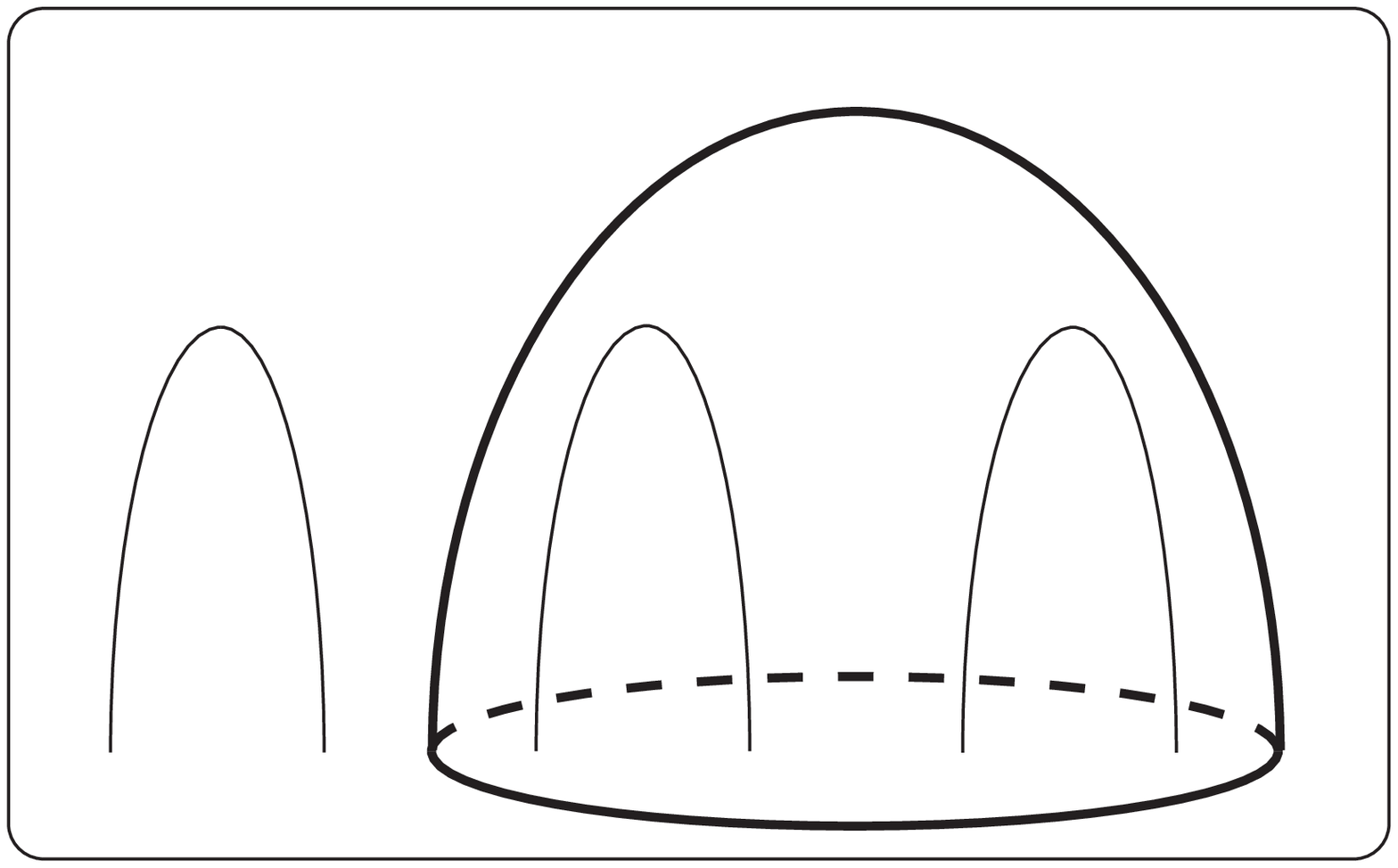} &
		\includegraphics[trim=0mm 0mm 0mm 0mm, width=.3\linewidth]{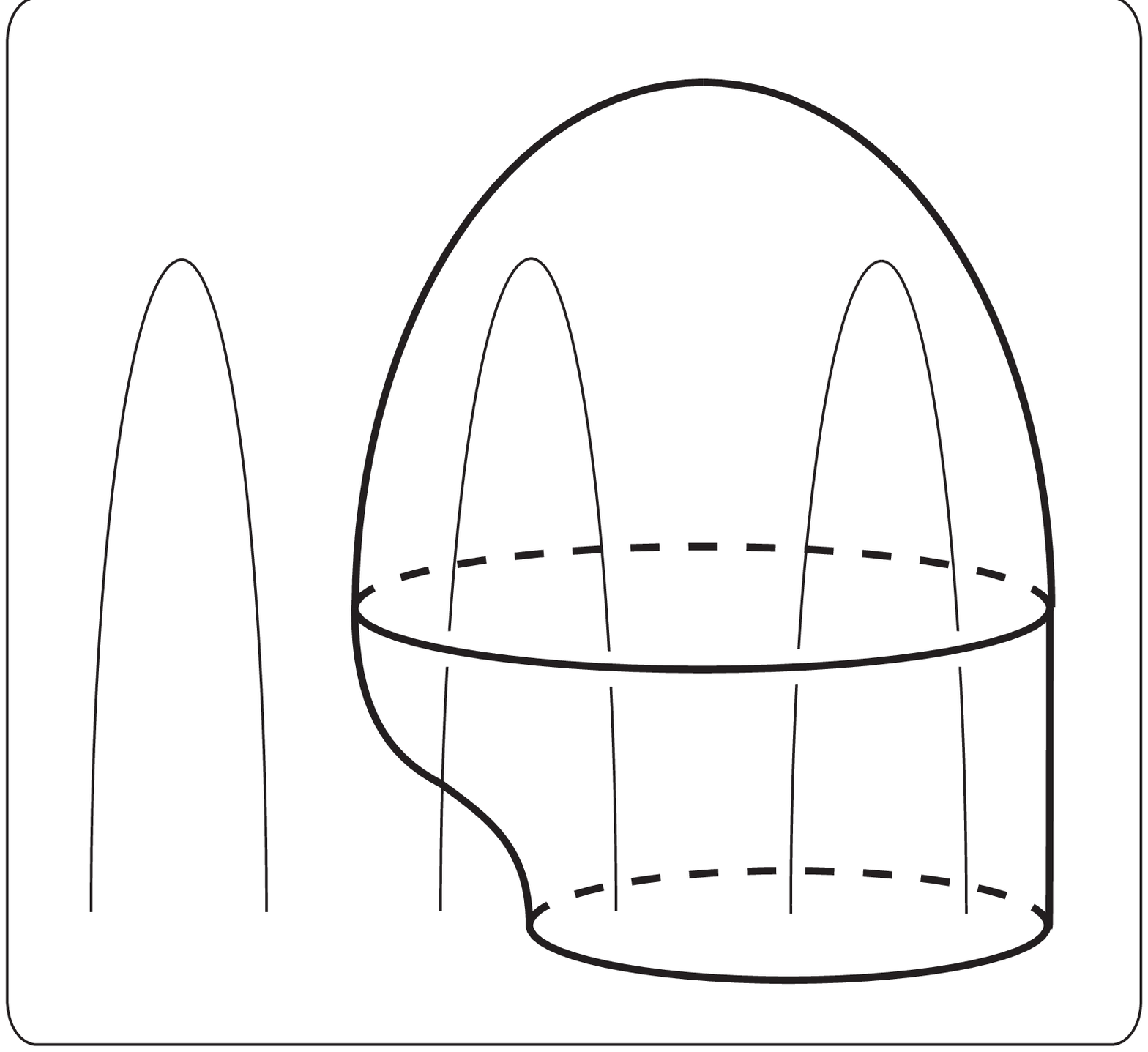} &
		\includegraphics[trim=0mm 0mm 0mm 0mm, width=.3\linewidth]{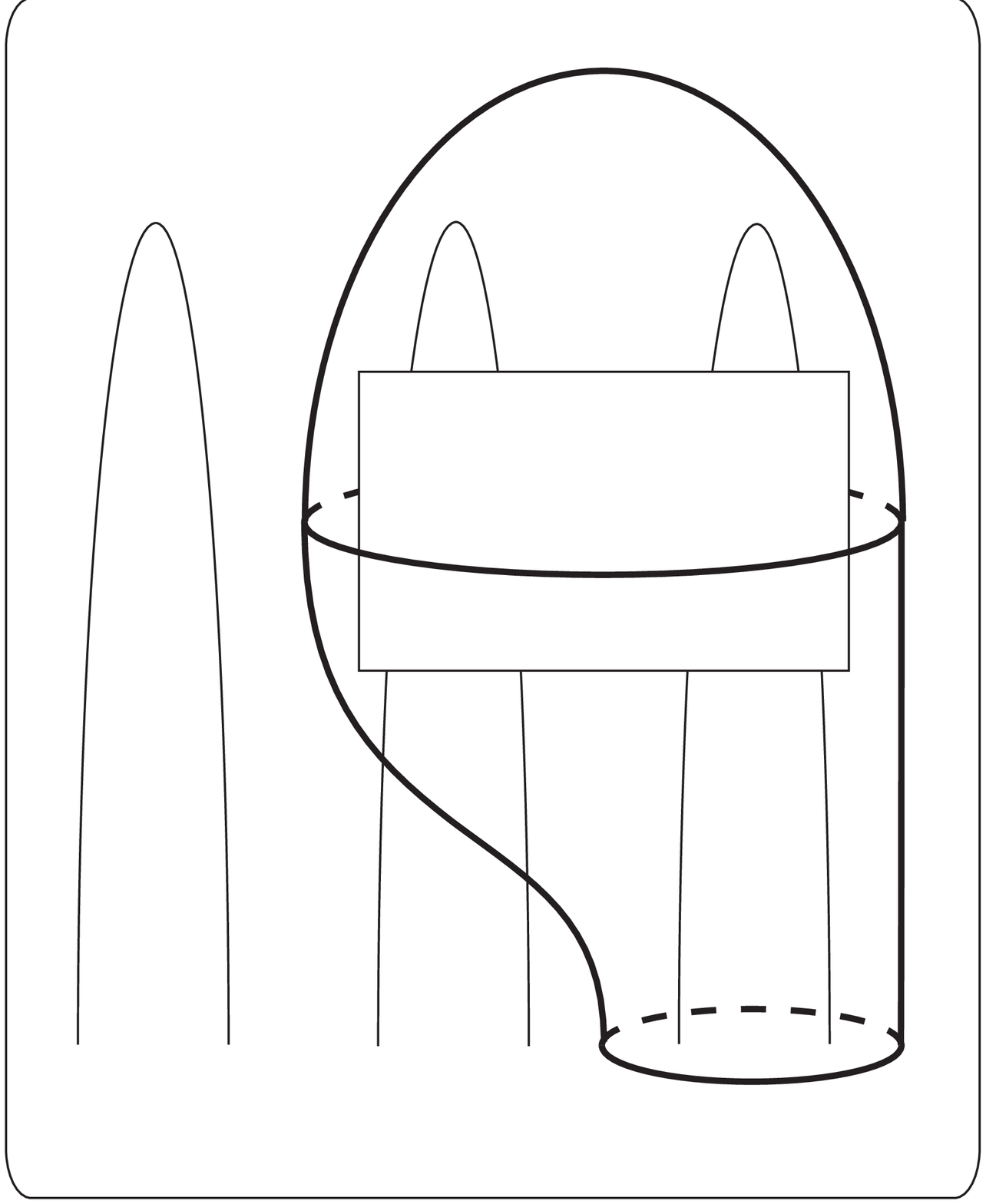}\\
		Type $D_0$ & Type $D_1$ & Type $D_2$
	\end{tabular}
	\end{center}
	\caption{}
	\label{D_0}
\end{figure}

For a 3-string trivial tangle $(B,T)$, we call the triple $(B,T,P)$ of the conclusion 1 in Lemma \ref{trivial} {\em Type $D_0$}, and the triple $(B,T,P)$ of the conclusion 2 {\em Type $D_1$}.
See Figure \ref{D_0}.

Next, we consider a disk $D$ intersecting $T$ in two points.
By Lemma \ref{trivial}, $D-T$ is $\partial$-compressible in $B-T$.
Hence, $D$ is obtained from two disks $D'$ and $D''$ each of which intersects $T$ in one point.
If both of $D'$ and $D''$ are of Type $D_1$, then $D-T$ is compressible in $B-T$.
If both of $D'$ and $D''$ are $\partial$-parallel in $(B,T)$, then $D-T$ is compressible in $B-T$ or $D$ is $\partial$-parallel in $(B,T)$.
Otherwise, $D$ is obtained from a disk of Type $D_1$ and a $\partial$-parallel disk.
We call the triple $(B,T,D)$ with such a disk $D$ {\em Type $D_2$}.
See Figure \ref{D_0}.
Here, the box indicates a 4-braid.

\begin{figure}[htbp]
	\begin{center}
		\includegraphics[trim=0mm 0mm 0mm 0mm, width=.5\linewidth]{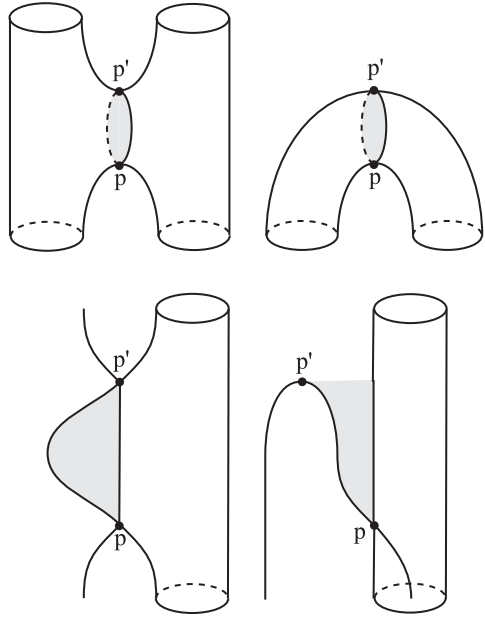}
	\end{center}
	\caption{gradient disks}
	\label{gradient}
\end{figure}

To describe the inside of the box exactly, we prepare the following term.
For a triple $(S^3,K,F)$ or $(B,T,P)$, let $h$ be a Morse function $h:S^3\to \Bbb{R}$ or $h:B\to \Bbb{R}$ such that $K|_h$ and $F|_h$, or $T|_h$ and $P|_h$ are also Morse functions.
Let $p$ be a saddle point of $F$ or $P$, or a point of intersection of $K\cap F$ or $T\cap P$.
Let $p'$ be a saddle or maximal point of $F$ or $P$, or a point of intersection of $K\cap F$ or $T\cap P$ or a maximal point of $K$ or $T$ that is immediately over $p$.
We say that a disk $D$ is an {\em (upper) gradient disk} for $p$ if $\rm{int}D|_h$ has no critical point, $\partial D-h^{-1}(h(p'))\subset F\cup K$ or $P\cup T$ and $\partial D$ passes through $p$ and $p'$.
See Figure \ref{gradient}.

\begin{figure}[htbp]
	\begin{center}
	\begin{tabular}{cc}
		\includegraphics[trim=0mm 0mm 0mm 0mm, width=.3\linewidth]{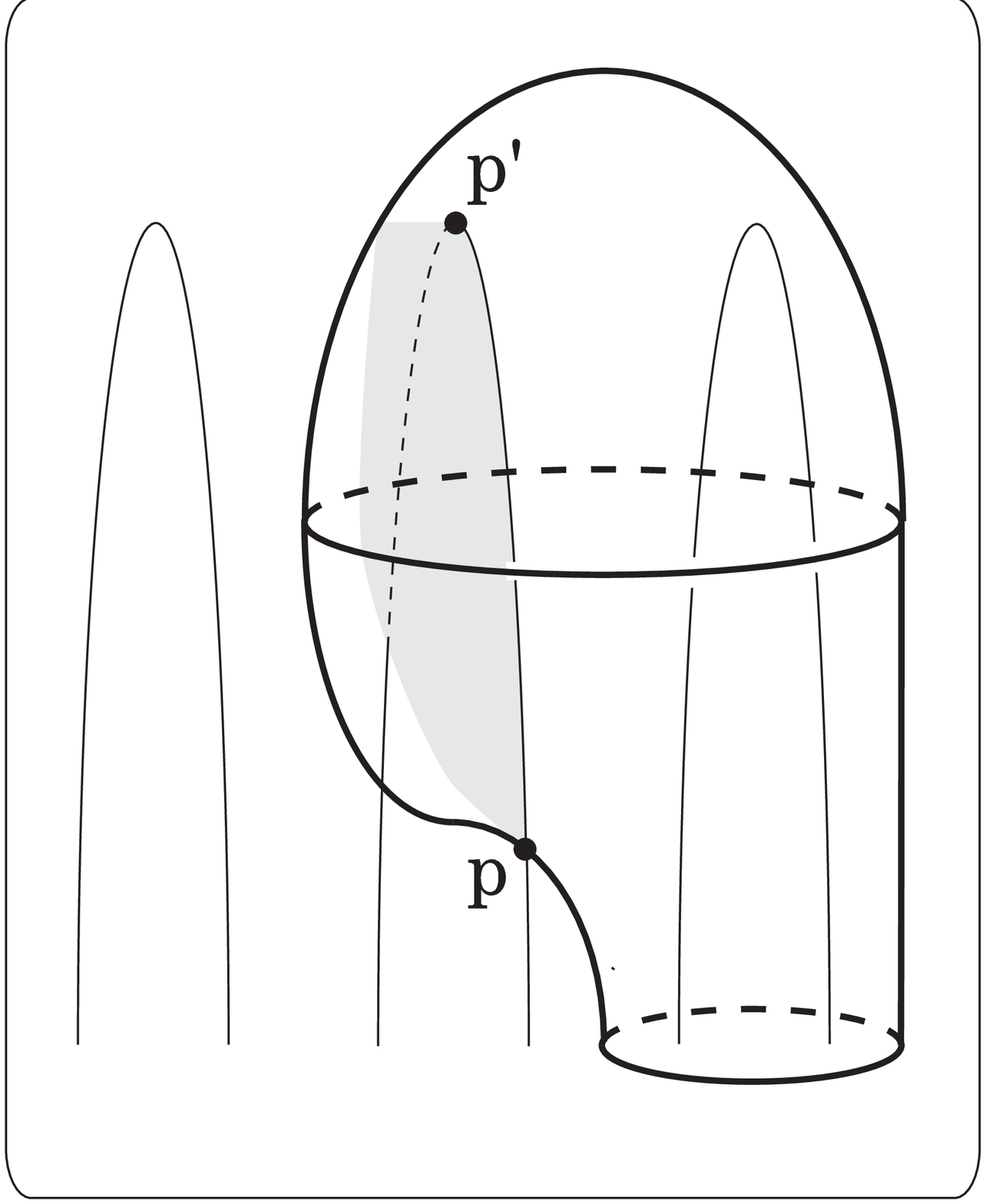} &
		\includegraphics[trim=0mm 0mm 0mm 0mm, width=.3\linewidth]{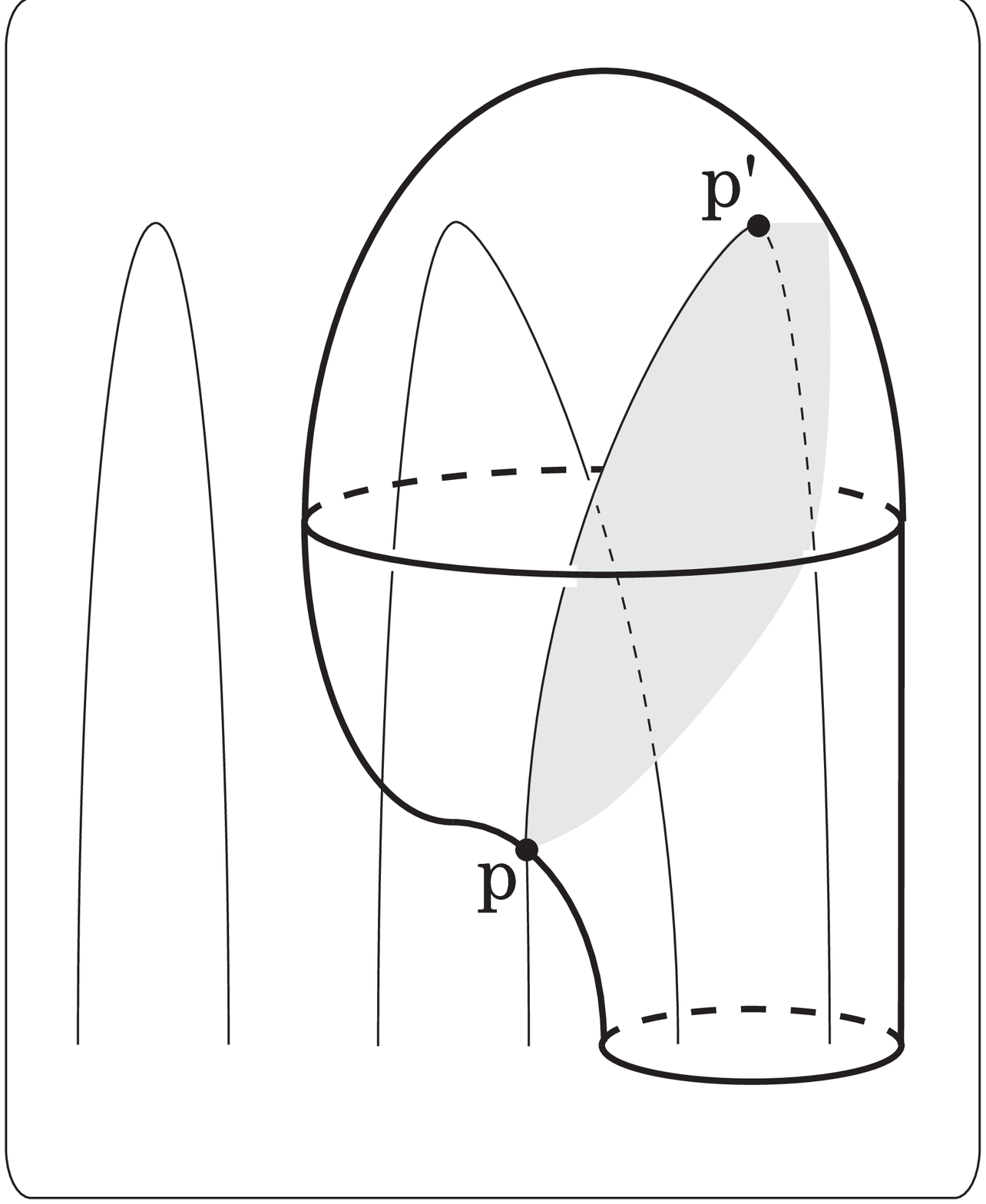} \\
		(a) & (b)
	\end{tabular}
	\end{center}
	\caption{gradient disks for Type $D_2$}
	\label{disk2_gradient}
\end{figure}

For the incompressibility of a disk $D$ of Type $D_2$, it is necessary that the 2-string subtangle including the 4-braid box is not split.
In other words, it can be express that there is no gradient disk passes through a point of $D\cap T$ and a maximal point of $T$ (Figre \ref{disk2_gradient}-(a)).
On the other side, if there exists another gradient disk passes through a point of $D\cap T$ and a maximal point of $T$, then $D$ is $\partial$-parallel in $(B,T)$ (Figre \ref{disk2_gradient}-(b)).
Here, we require that two maximal points of $T$ are in a same level.
In summary, a disk $D$ of Type $D_2$ is incompressible, meridionally incompressible and not $\partial$-parallel in $(B,T)$ if and only if there exists no gradient disk for two points of $D\cap T$.

Same arguments as above show that;

\begin{lemma}
Let $(B,T)$ be a 3-string trivial tangle and $P$ be an incompressible, meridionally incompressible and not $\partial$-parallel surface in $(B,T)$.
If $P$ is;
\begin{enumerate}
	\item a disk which intersects $T$ in $k$ points, then the triple $(B,T,P)$ is Type $D_k$ $(k=0,1,2,3)$.
	\item an annulus which intersects $T$ in $k$ points, then the triple $(B,T,P)$ is Type $A_k$ $(k=0,1)$.
	\item an annulus which intersects $T$ in 2 points, then the triple $(B,T,P)$ is Type $A_{21}$ or $A_{22}$ or $A_{23}$.
	\item a pants which is disjoint from $T$, then the triple $(B,T,P)$ is Type $P_0$.
	\item a pants which intersects $T$ in 1 point, then the triple $(B,T,P)$ is Type $P_{11}$ or $P_{12}$.
	\item a 4-punctured sphere which is disjoint from $T$, then the triple $(B,T,P)$ is Type $Q_{01}$ or $Q_{02}$.
	\item a once punctured torus which is disjoint from $T$, then the triple $(B,T,P)$ is Type $T_k$ $(k=0,1)$.
	\item a twice punctured torus which is disjoint from $T$, then the triple $(B,T,P)$ is Type $U_{02}$ or $U_{03}$.
\end{enumerate}
\end{lemma}

\bigskip

\begin{tabular}{c|cccc}
$P$ $\backslash $ $|P\cap T|$ & 0 & 1 & 2 & 3\\
\hline
disk & $D_0$ & $D_1$ & $D_2$ & $D_3$\\
annulus & $A_0$ & $A_1$ & $A_{21}$, $A_{22}$, $A_{23}$ & \\
pants & $P_0$ & $P_{11}$, $P_{12}$ & & \\
4-punctured sphere & $Q_{01}$, $Q_{02}$ & & & \\
once punctured torus & $T_0$ & $T_1$ & & \\
twice punctured torus & $U_{02}$, $U_{03}$ & & 
\end{tabular}
\smallskip
\begin{center}
Table of Types
\end{center}



\begin{figure}[htbp]
\begin{center}
\begin{tabular}{ccc}
	\includegraphics[trim=0mm 0mm 0mm 0mm, width=.3\linewidth]{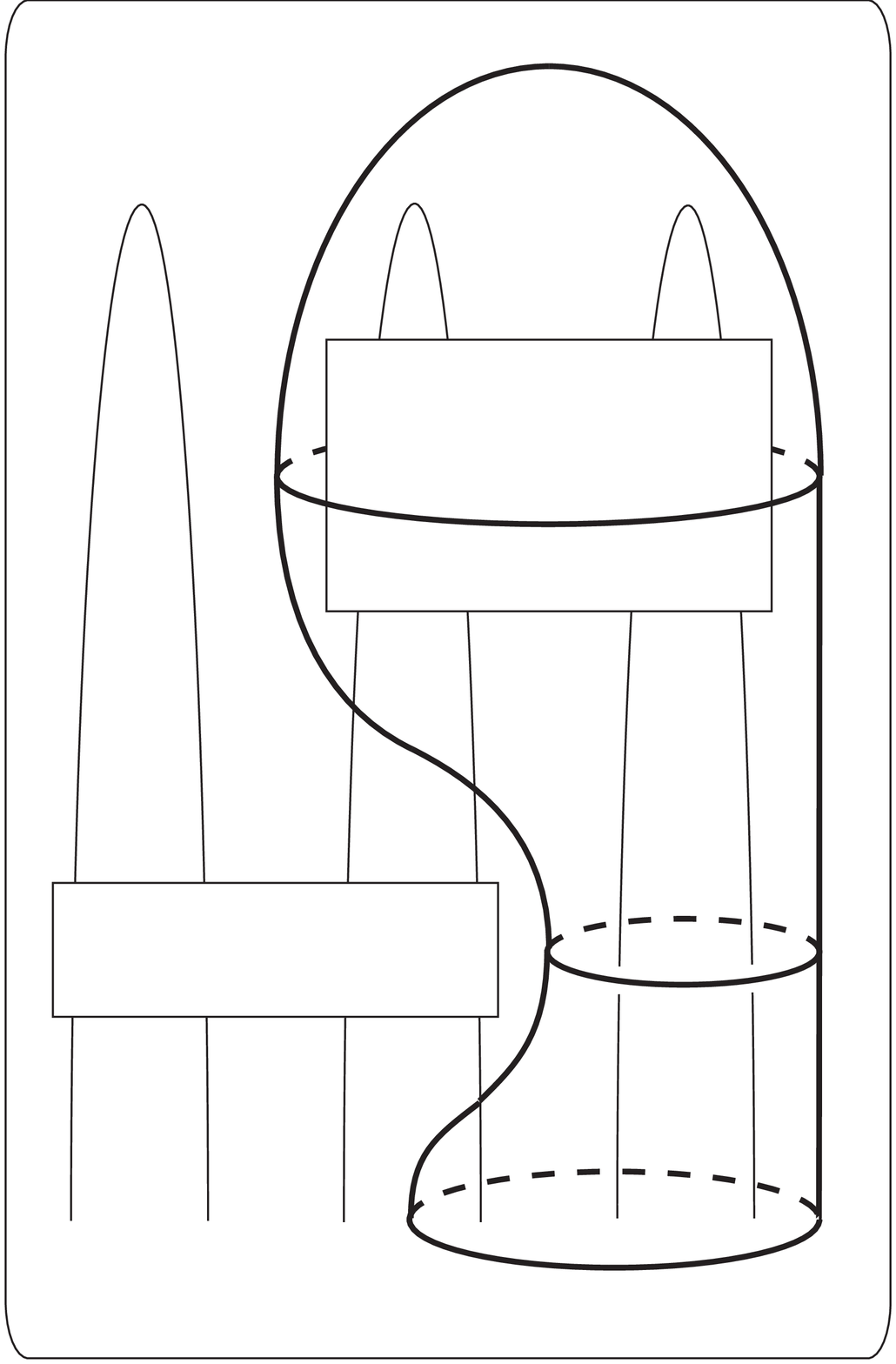}&
	\includegraphics[trim=0mm 0mm 0mm 0mm, width=.3\linewidth]{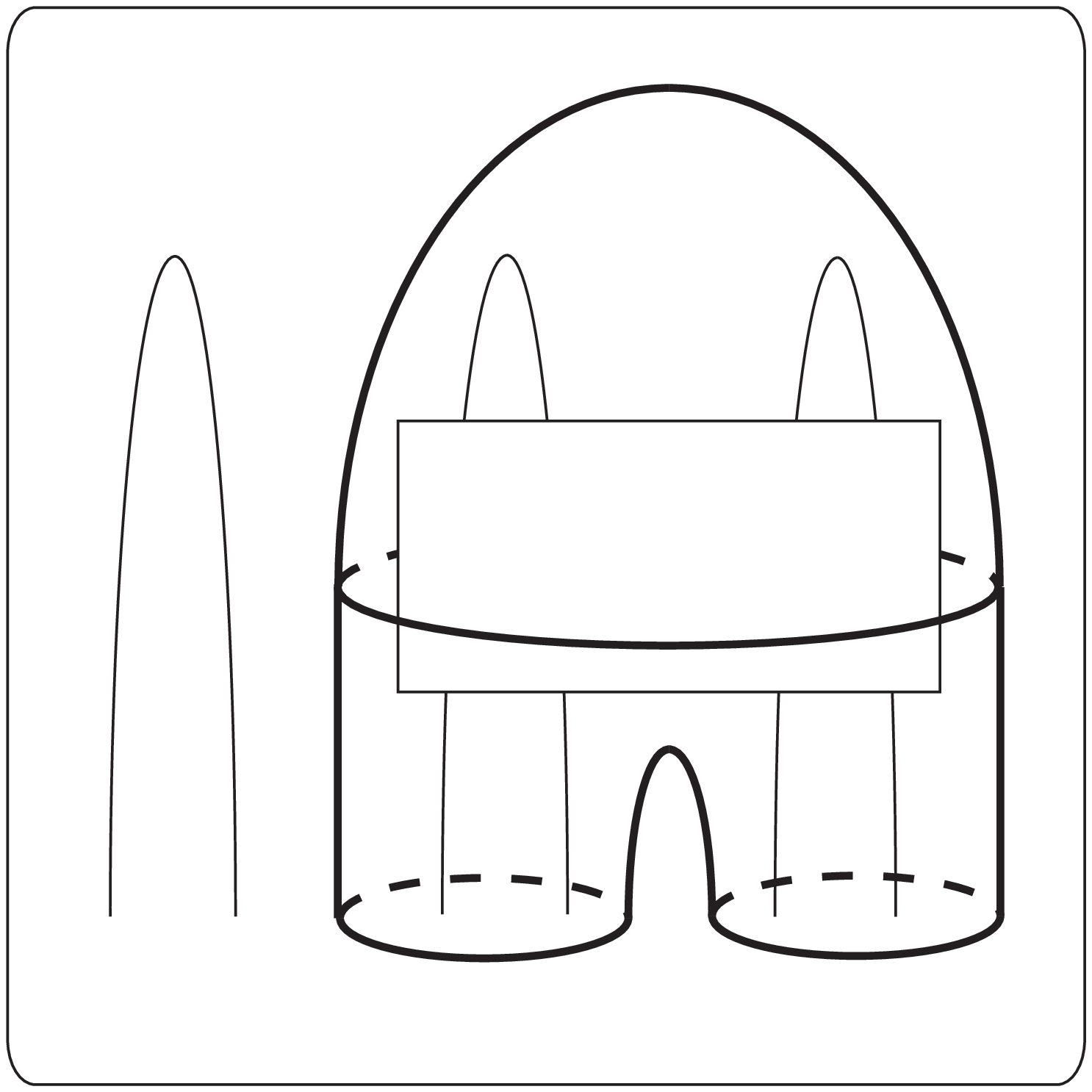}&
	\includegraphics[trim=0mm 0mm 0mm 0mm, width=.3\linewidth]{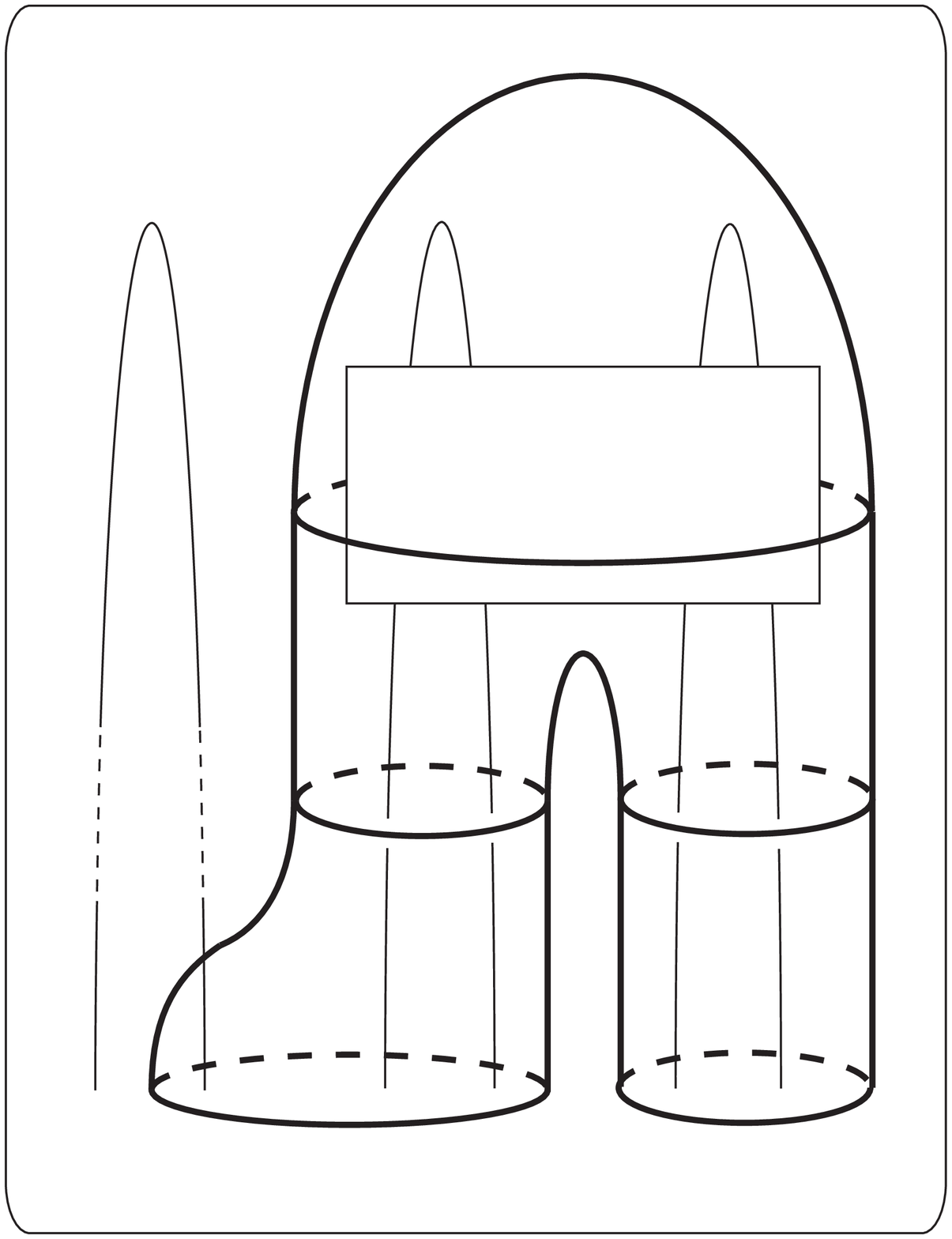}\\
	Type $D_3$ & Type $A_0$ & Type $A_1$
	\end{tabular}
\end{center}
	\caption{}
	\label{}
\end{figure}

\begin{figure}[htbp]
\begin{center}
\begin{tabular}{ccc}
	\includegraphics[trim=0mm 0mm 0mm 0mm, width=.3\linewidth]{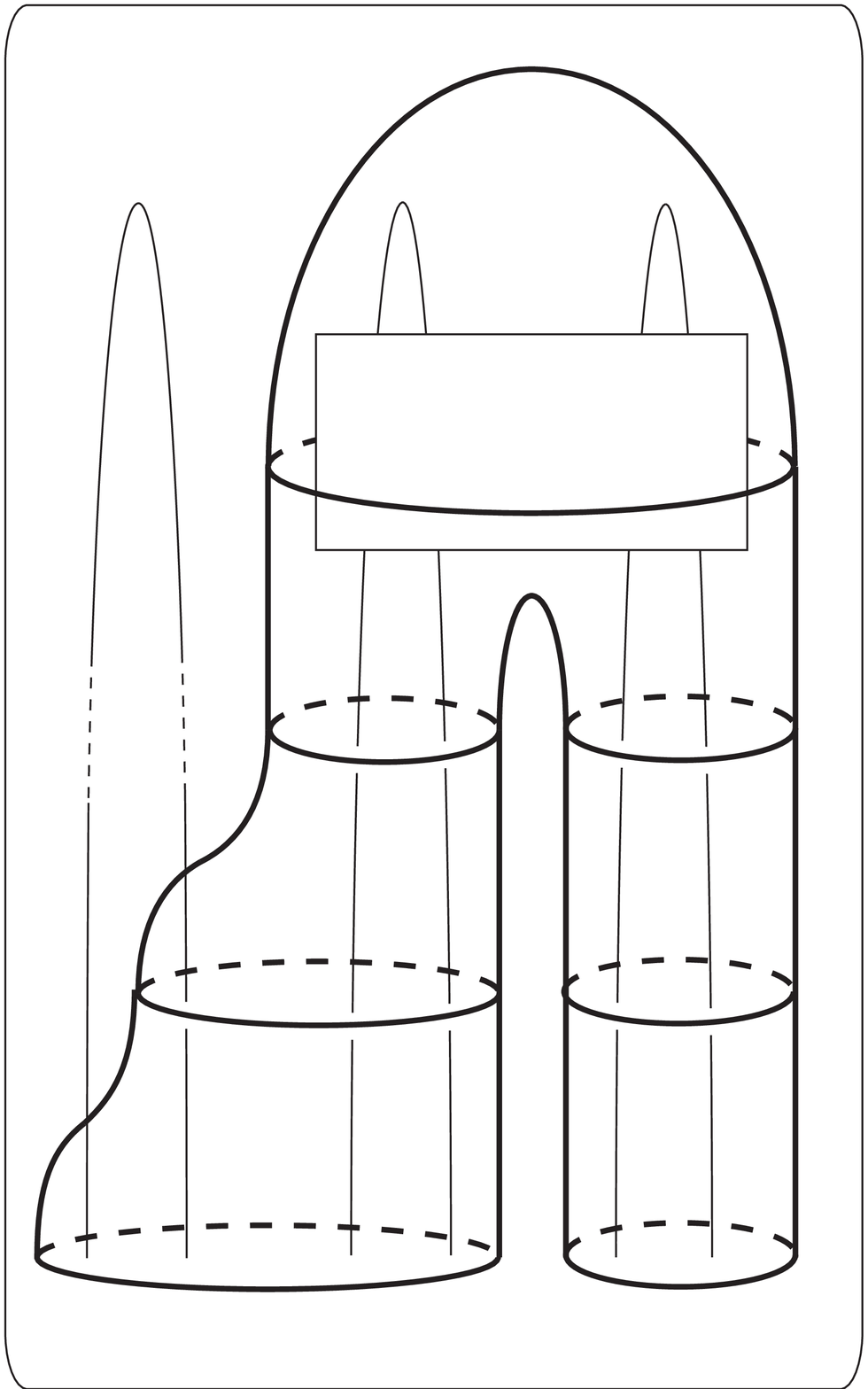}&
	\includegraphics[trim=0mm 0mm 0mm 0mm, width=.3\linewidth]{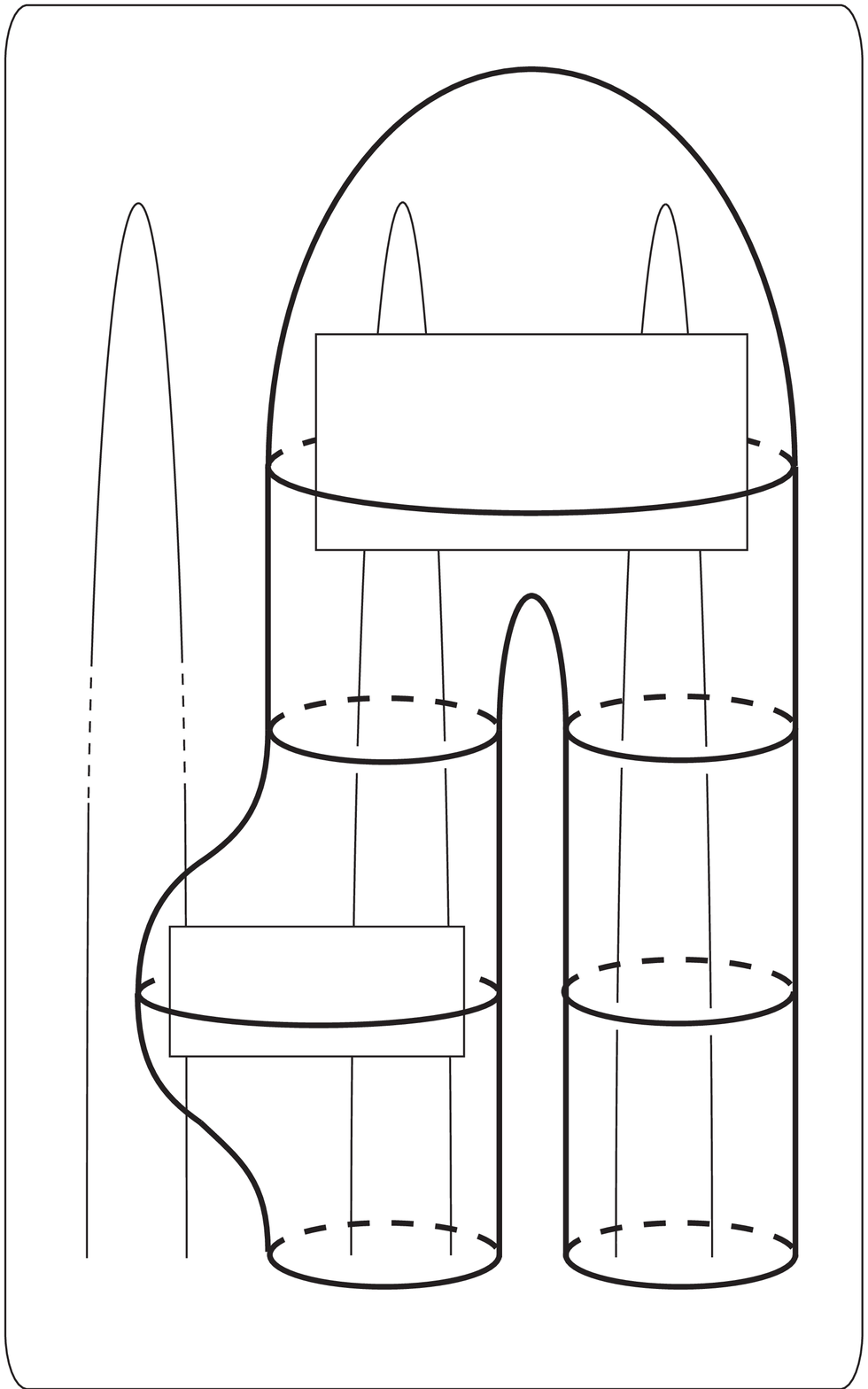}&
	\includegraphics[trim=0mm 0mm 0mm 0mm, width=.3\linewidth]{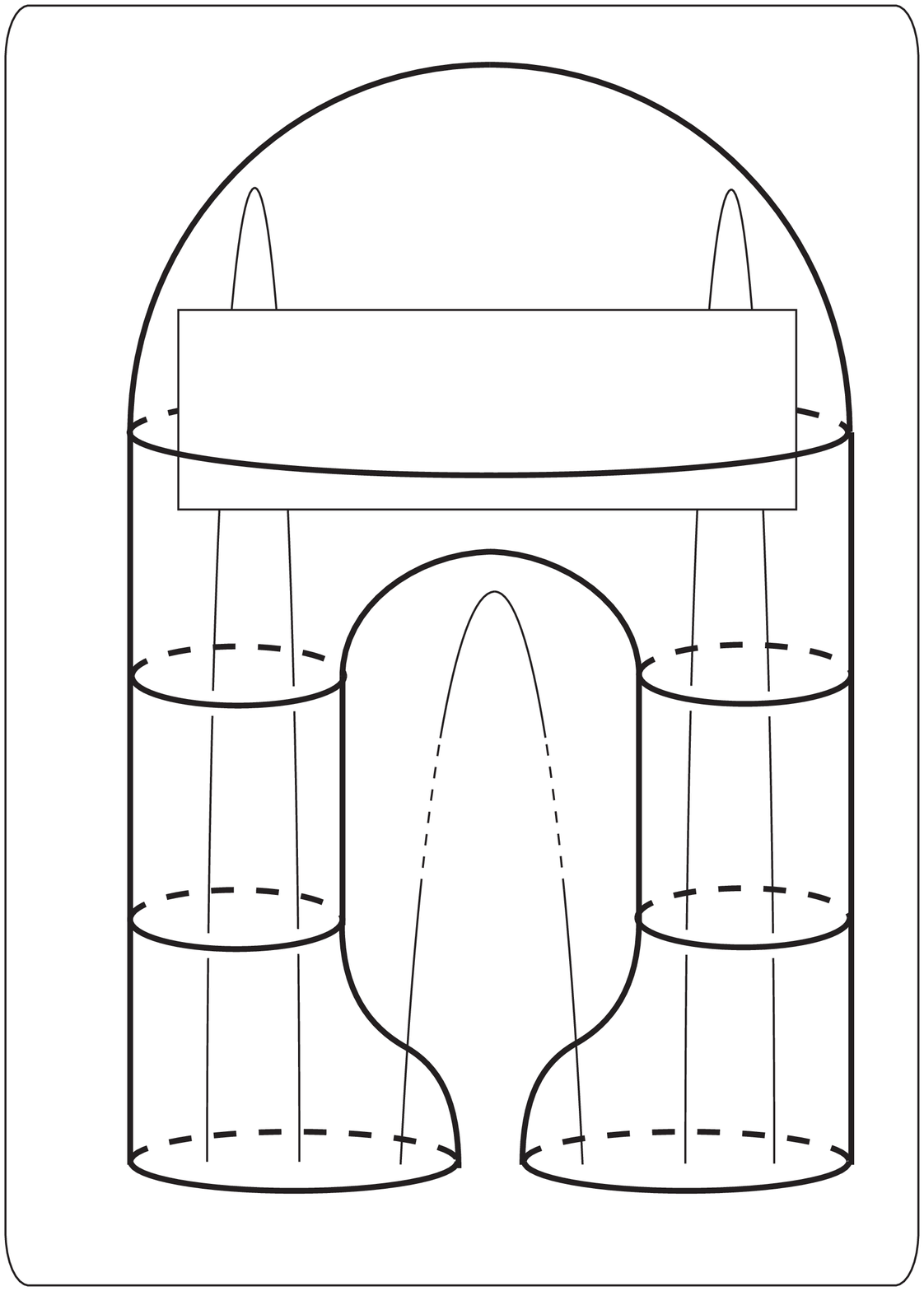}\\
	nonexistence & Type $A_{21}$ & Type $A_{22}$
	\end{tabular}
\end{center}
	\caption{}
	\label{non}
\end{figure}

\begin{remark}
In Figure \ref{non}, the left-hand annulus can not be incompressible.
The string in the outside of the annulus can be untied by sliding on ``instep'' of the annulus.
\end{remark}

\begin{figure}[htbp]
\begin{center}
\begin{tabular}{ccc}
	\includegraphics[trim=0mm 0mm 0mm 0mm, width=.3\linewidth]{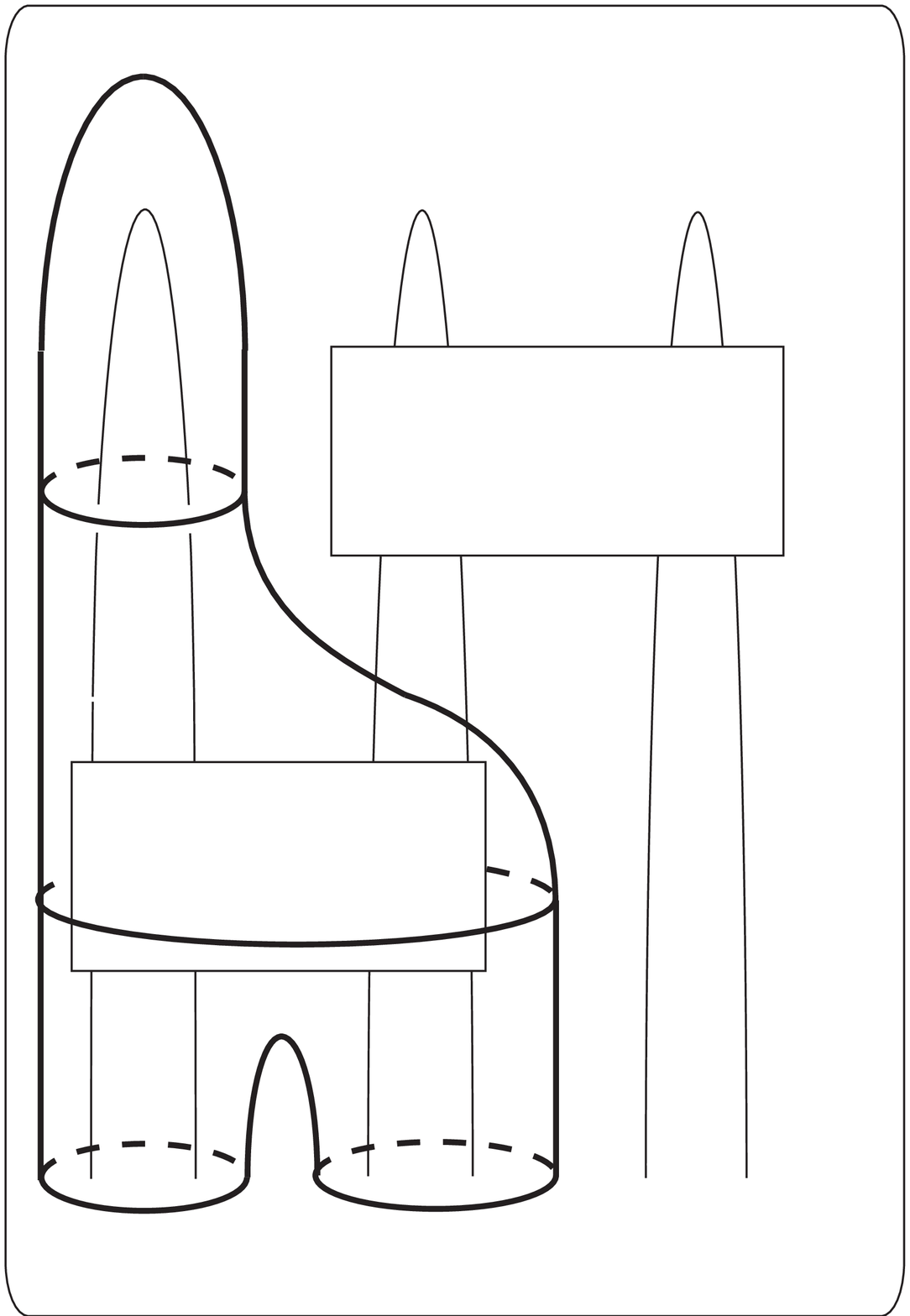}&
	\includegraphics[trim=0mm 0mm 0mm 0mm, width=.3\linewidth]{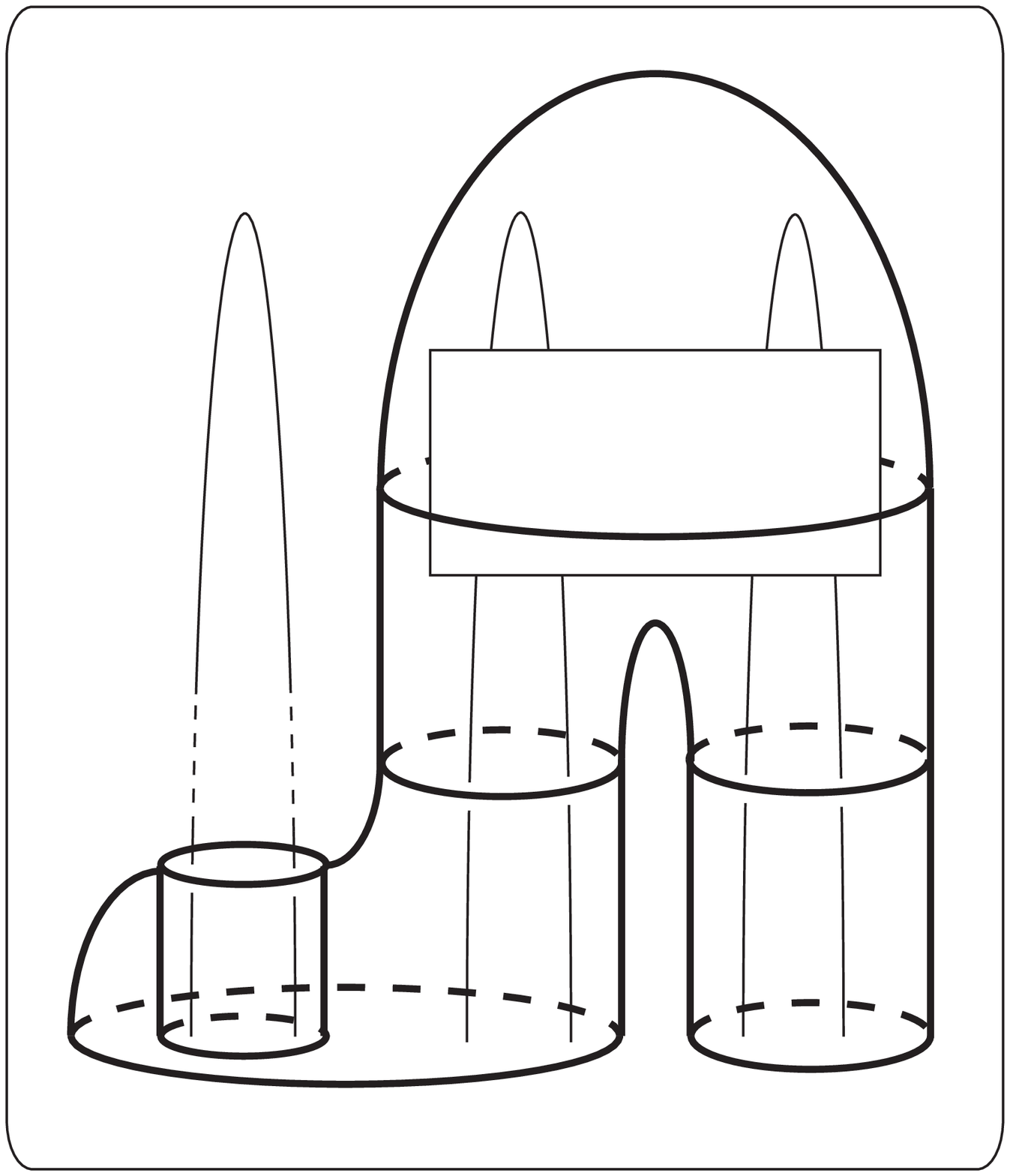}&
	\includegraphics[trim=0mm 0mm 0mm 0mm, width=.3\linewidth]{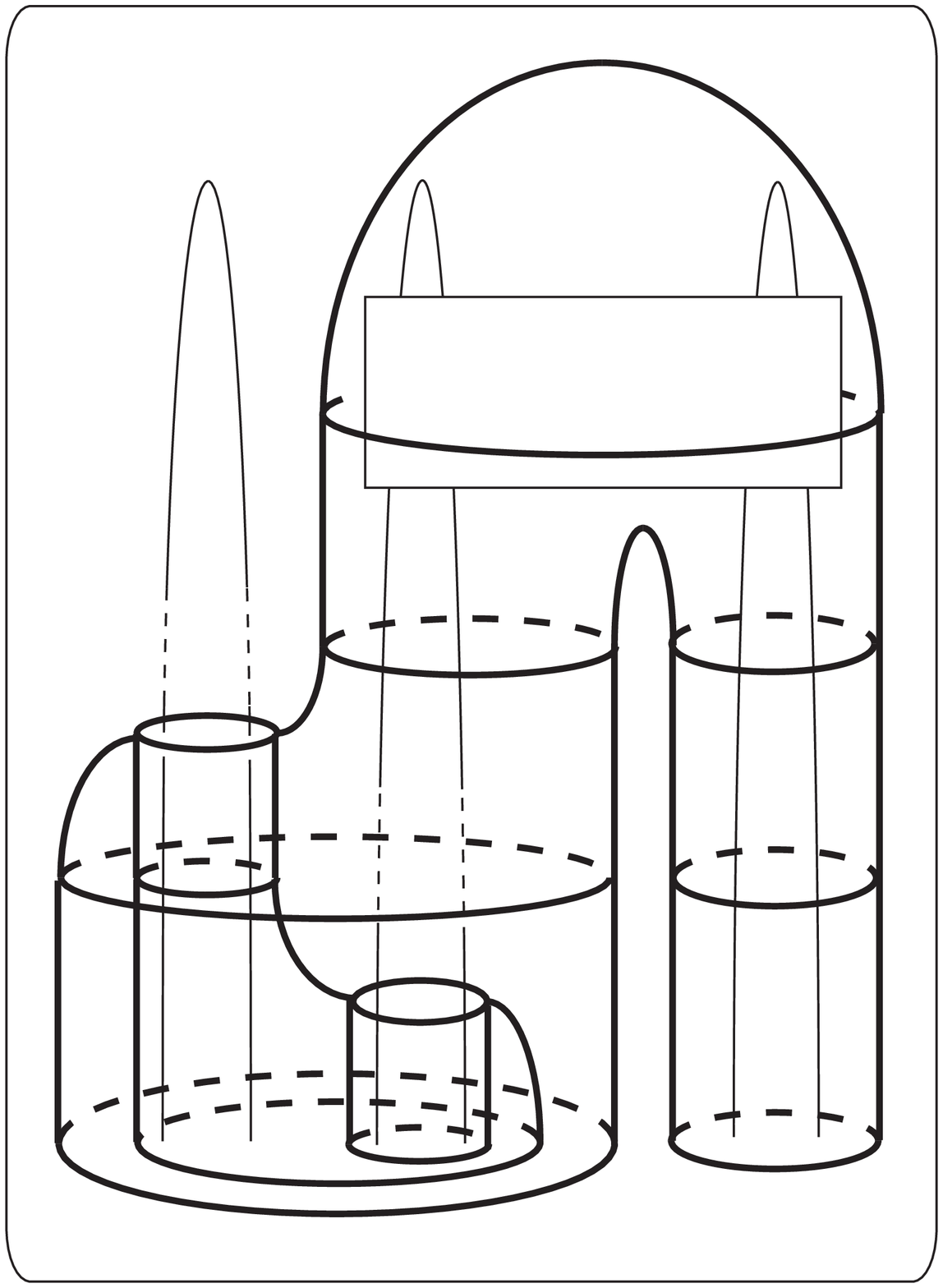}\\
	Type $A_{23}$ & Type $P_0$ & Type $Q_{01}$
	\end{tabular}
\end{center}
	\caption{}
	\label{}
\end{figure}


\begin{figure}[htbp]
\begin{center}
\begin{tabular}{ccc}
	\includegraphics[trim=0mm 0mm 0mm 0mm, width=.3\linewidth]{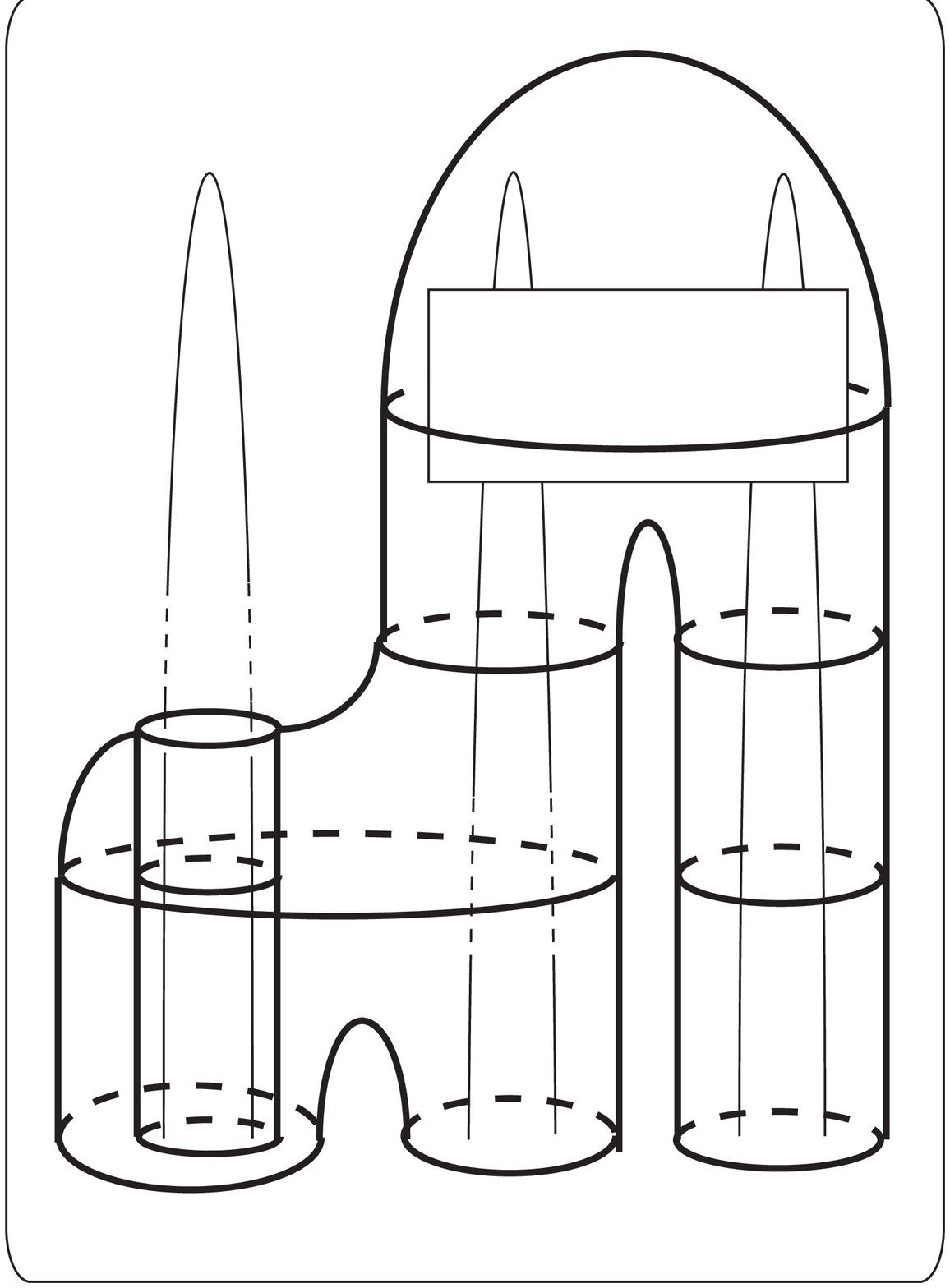}&
	\includegraphics[trim=0mm 0mm 0mm 0mm, width=.3\linewidth]{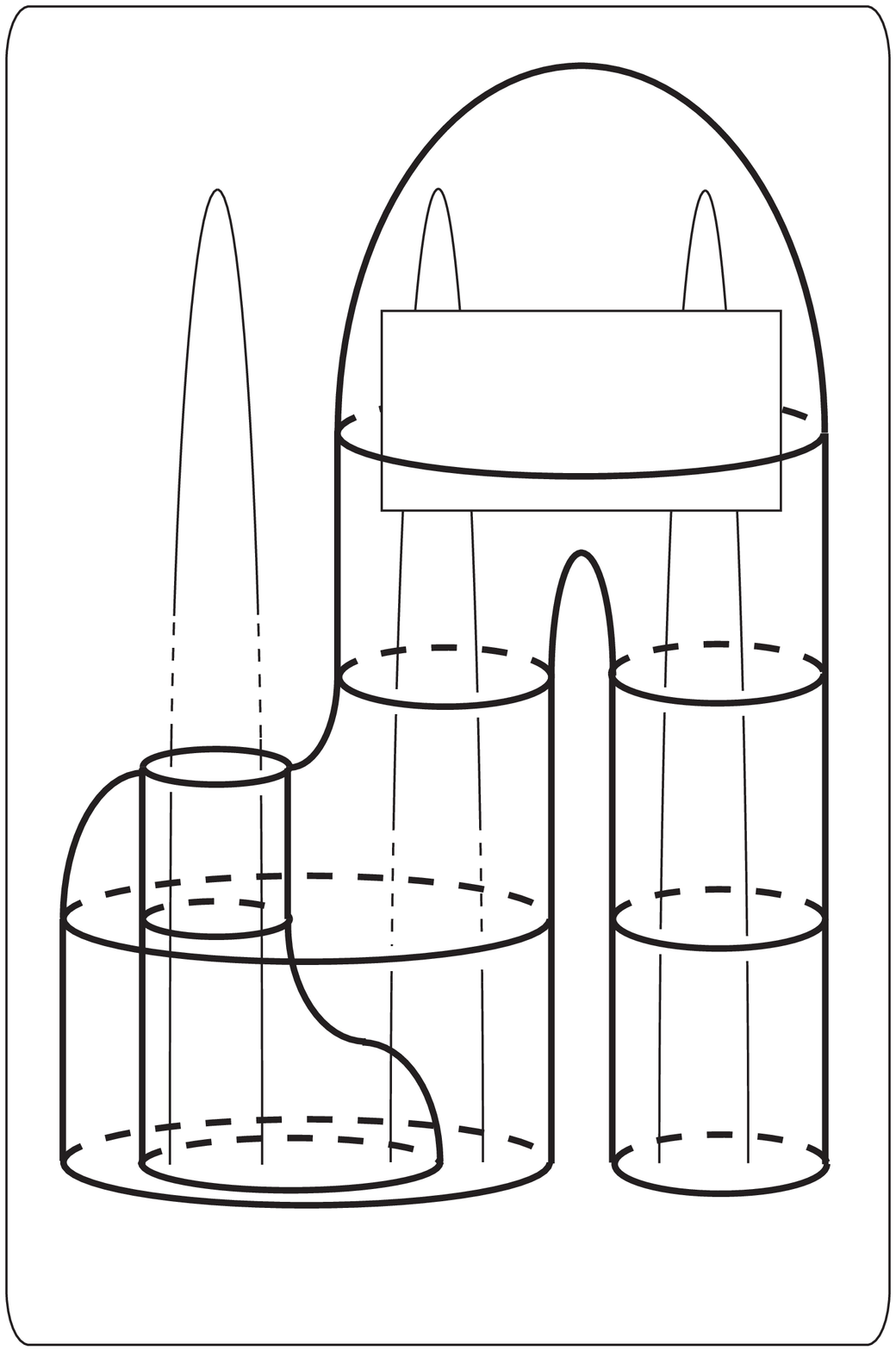}&
	\includegraphics[trim=0mm 0mm 0mm 0mm, width=.3\linewidth]{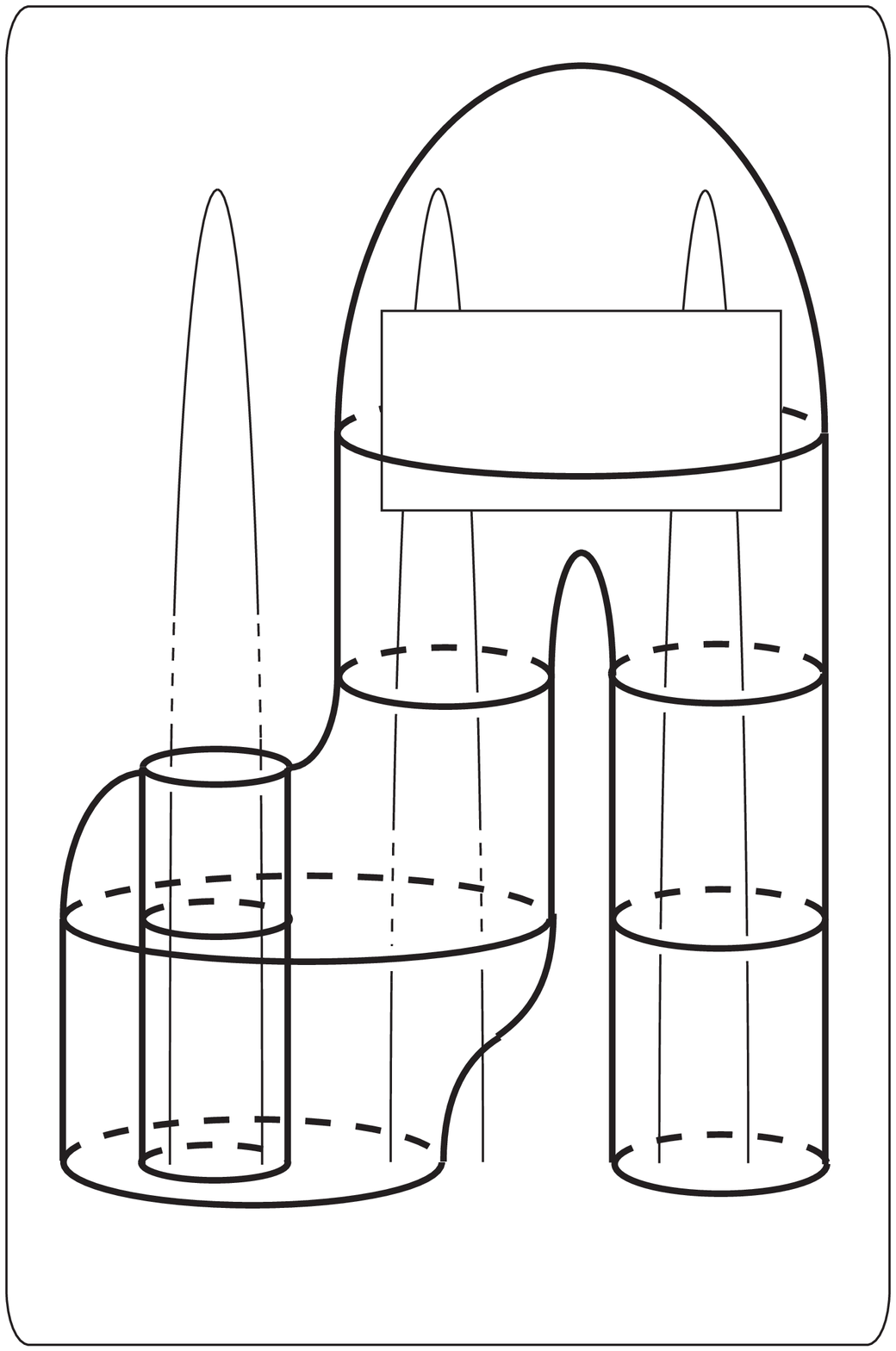}\\
	Type $Q_{02}$ & Type $P_{11}$ & Type $P_{12}$
	\end{tabular}
\end{center}
	\caption{}
	\label{}
\end{figure}

\begin{figure}[htbp]
\begin{center}
\begin{tabular}{ccc}
	\includegraphics[trim=0mm 0mm 0mm 0mm, width=.3\linewidth]{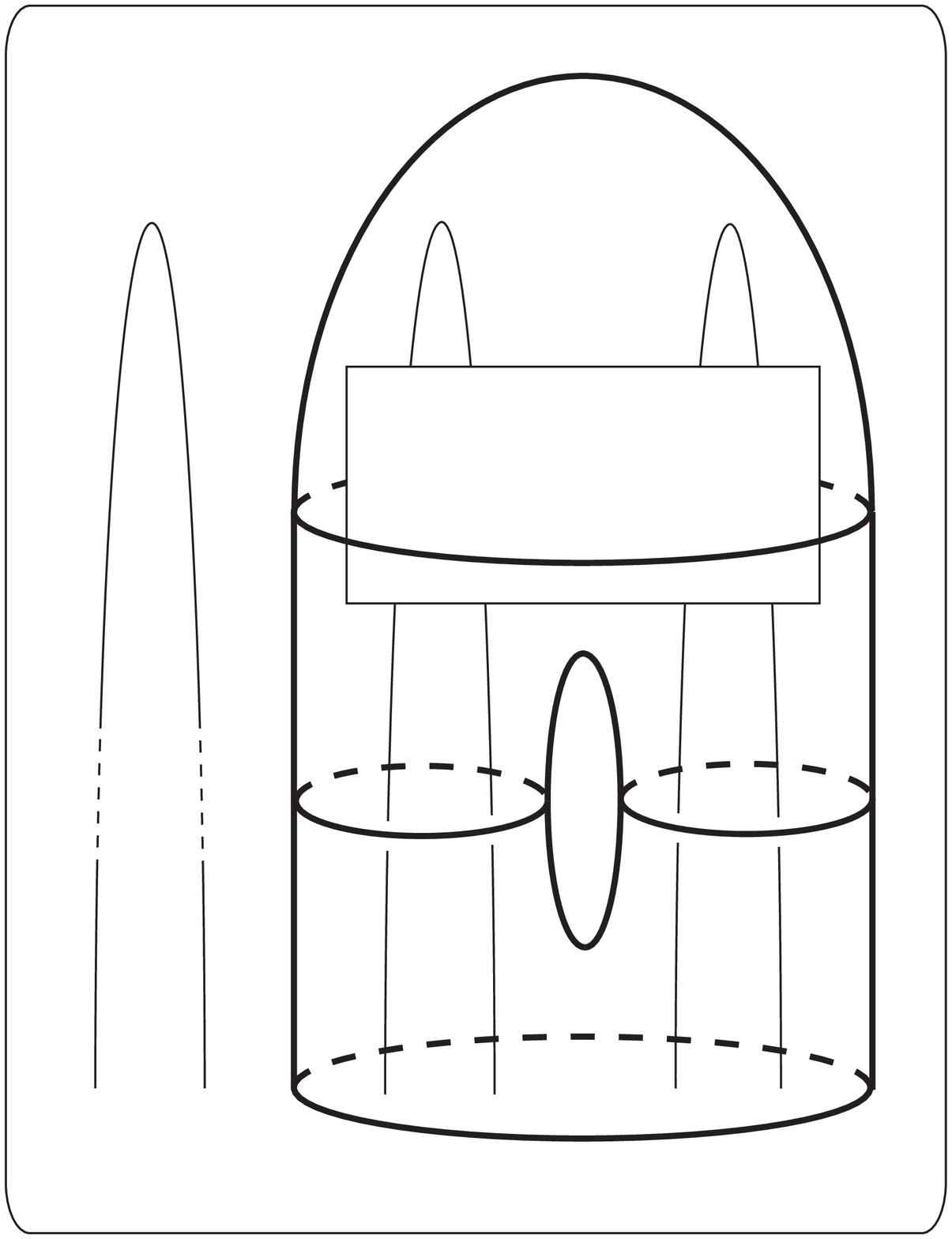}&
	\includegraphics[trim=0mm 0mm 0mm 0mm, width=.3\linewidth]{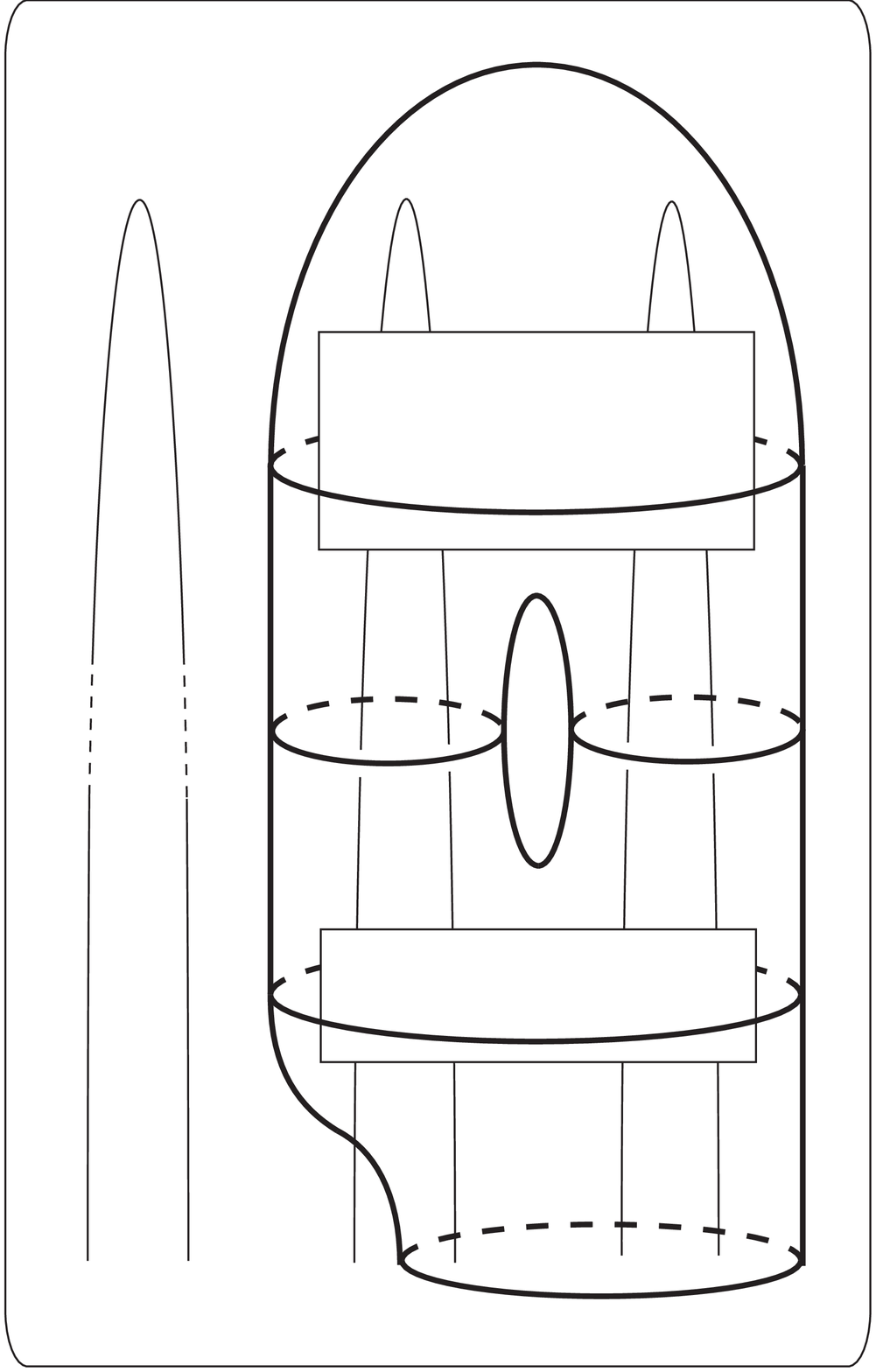}&
	\includegraphics[trim=0mm 0mm 0mm 0mm, width=.3\linewidth]{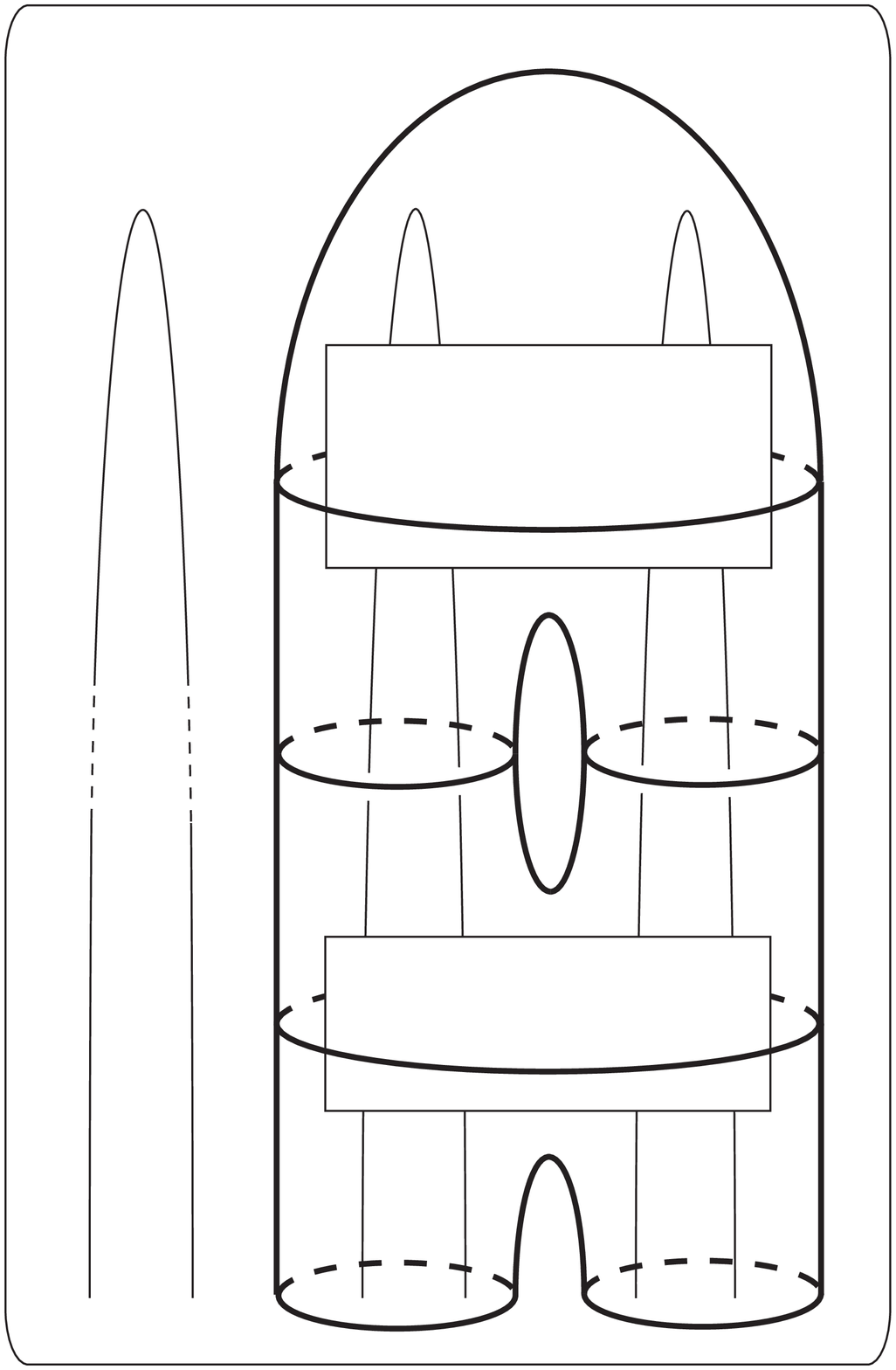}\\
	Type $T_{0}$ & Type $T_1$ & Type $U_{02}$
	\end{tabular}
\end{center}
	\caption{}
	\label{}
\end{figure}

\begin{figure}[htbp]
\begin{center}
\begin{tabular}{cc}
	\includegraphics[trim=0mm 0mm 0mm 0mm, width=.3\linewidth]{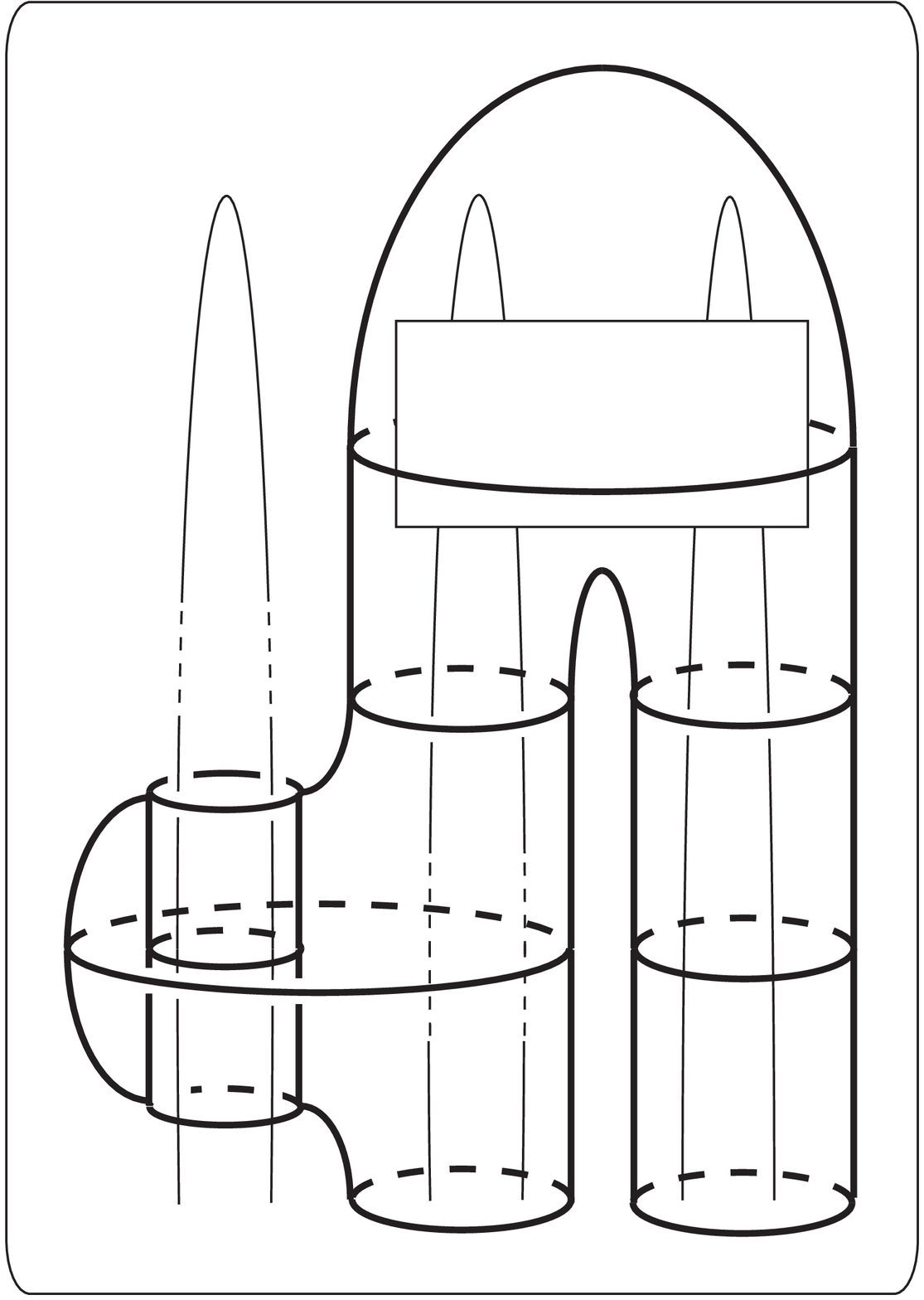}&
	\includegraphics[trim=0mm 0mm 0mm 0mm, width=.3\linewidth]{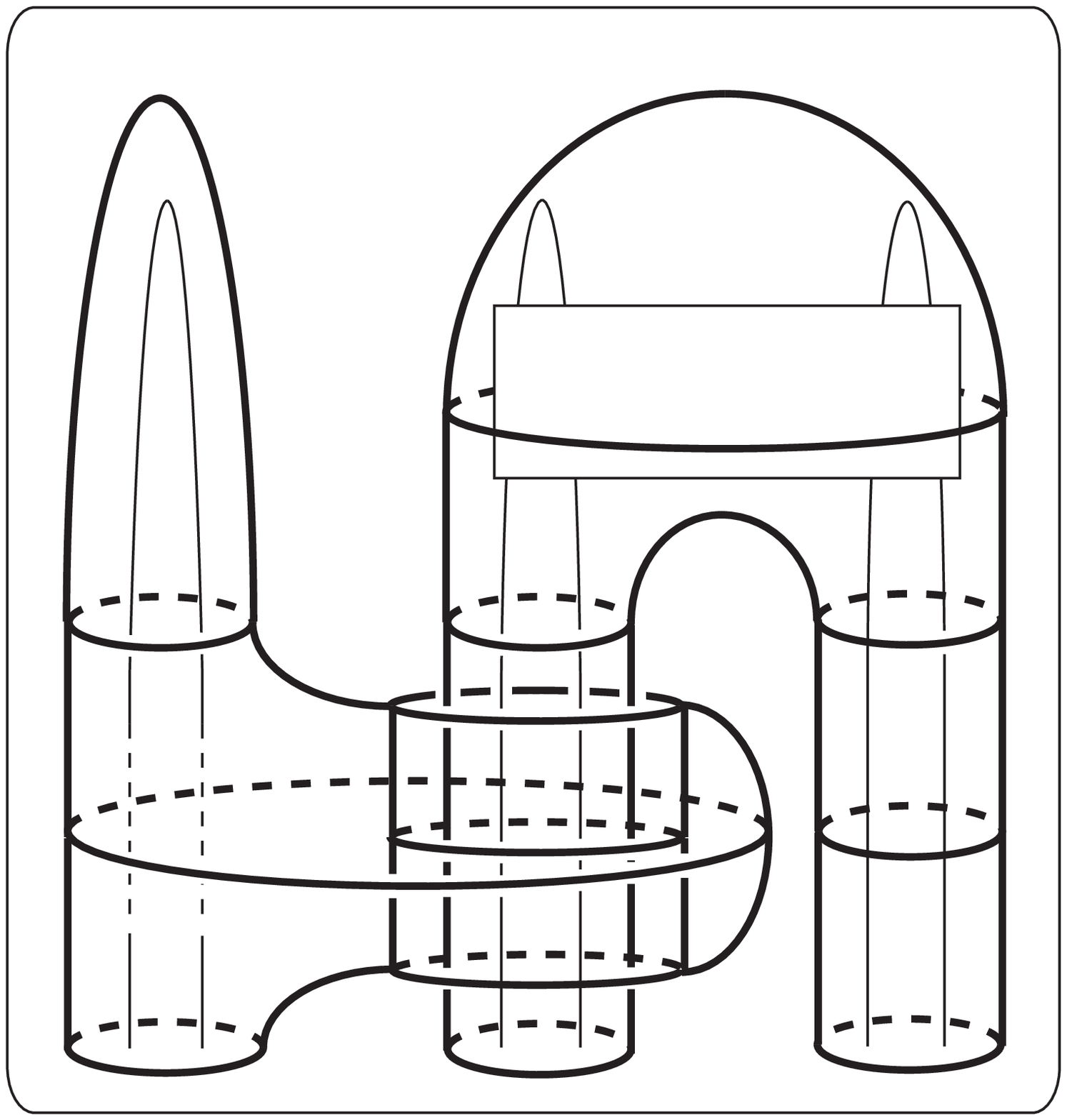}\\
	Type $U_{03}$ & Type $A_0\cup T_0$
	\end{tabular}
\end{center}
	\caption{}
	\label{U_{03}}
\end{figure}

\begin{figure}[htbp]
\begin{center}
	\includegraphics[trim=0mm 0mm 0mm 0mm, width=.7\linewidth]{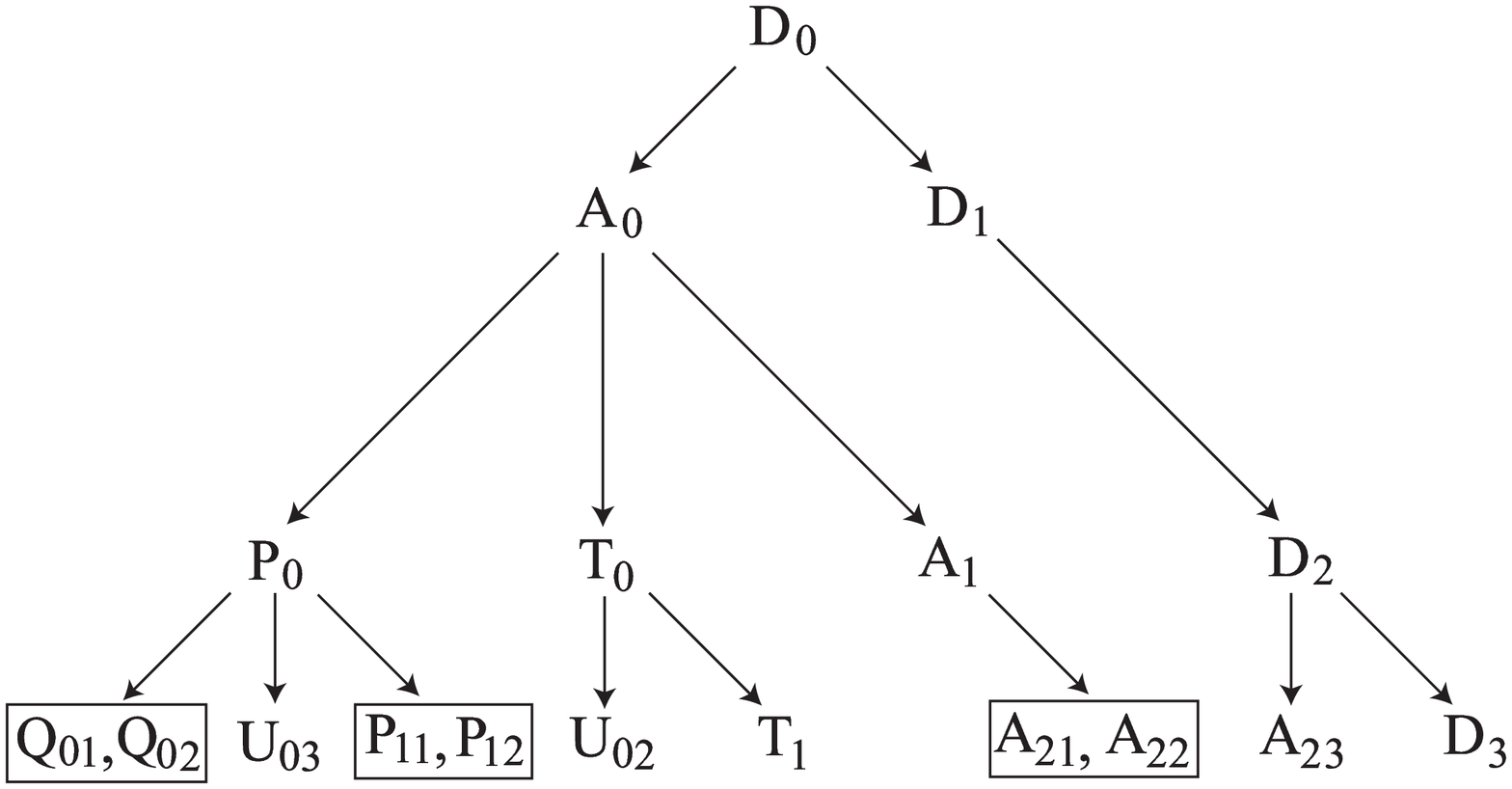}
\end{center}
	\caption{A family tree}
	\label{family}
\end{figure}

By Lemma \ref{trivial}, any incompressible and meridionally incompressible surface in a 3-string trivial tangle, except for a disk of Type $D_0$ or $D_1$, has a $\partial$-compressing disk, and hence has a ``parent'', that is, a surface obtained by a $\partial$-compression.
Figure \ref{family} expresses a family tree for these types.

Every two surfaces of these types may not coexist in a 3-string trivial tangle.
For example, Type $A_0$ and $U_{02}$ cannot coexist.
Otherwise, the bottom braid box for a surface of Type $U_{02}$ must be unknotted and it becomes to be compressible.
When there are two types $X$ and $Y$, we denote the type by {\em Type $X\cup Y$}.
For example, if there are an annlus of Type $A_0$ and a torus of $T_0$, then we have Type $A_0\cup T_0$ as Figure \ref{U_{03}}.

Here, we enumerate all pairs of two Types that can coexist.

\begin{lemma}\label{coexist}
Let $X$ and $Y$ be one of $D_0$, $D_1$, $D_2$, $D_3$, $A_0$, $A_1$, $A_{21}$, $A_{22}$, $A_{23}$, $P_0$, $P_{11}$, $P_{12}$, $Q_{01}$, $Q_{02}$, $T_0$, $T_1$, $U_{02}$, $U_{03}$.
Then, Type $X$ and Type $Y$ can coexist if and only if $X\cup Y$ is either ;
$D_0\cup D_1$, $D_0\cup D_2$, $D_0\cup A_0$, $D_1\cup D_2$, $D_2\cup D_3$, $D_2\cup A_0$, $D_2\cup A_1$, $D_2\cup A_{21}$, $D_2\cup A_{22}$, $D_2\cup A_{23}$, $D_2\cup P_0$, $D_2\cup Q_{01}$, $D_2\cup Q_{02}$, $D_2\cup P_{11}$, $D_2\cup P_{12}$, $D_2\cup T_{0}$, $D_2\cup T_{03}$, $D_3\cup A_{1}$, $D_3\cup A_{22}$, $D_3\cup T_{0}$, $A_0\cup A_{1}$, $A_0\cup A_{22}$, $A_0\cup P_{0}$, $A_0\cup T_{0}$, $A_1\cup A_{21}$, $A_1\cup A_{22}$, $A_{21}\cup A_{22}$, $P_0\cup T_0$, $P_0\cup Q_{01}$, $P_0\cup Q_{02}$, $P_0\cup P_{11}$, $P_0\cup P_{12}$, $P_0\cup T_{03}$, $Q_{01}\cup T_{03}$, $Q_{02}\cup T_{03}$, $T_0\cup P_{11}$, $T_0\cup P_{12}$, $T_0\cup T_{1}$, $T_0\cup T_{02}$.
\end{lemma}

\begin{figure}[htbp]
\begin{center}
	\includegraphics[trim=0mm 0mm 0mm 0mm, width=.8\linewidth]{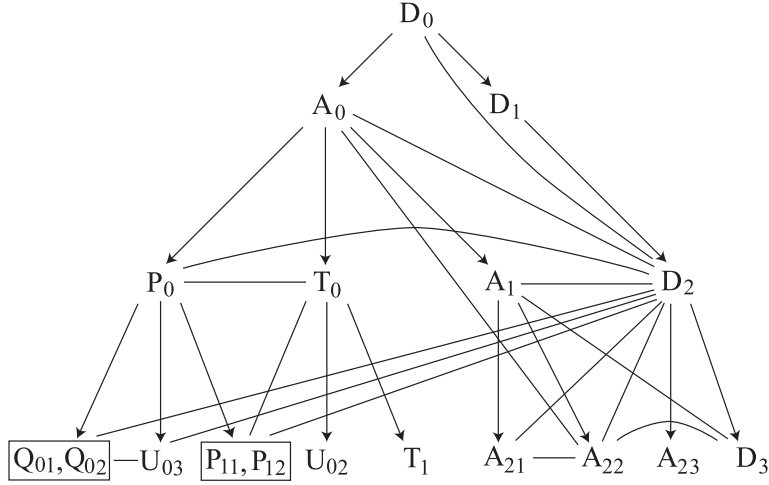}
\end{center}
	\caption{coexistable pairs}
	\label{coexist_fig}
\end{figure}

Figure \ref{coexist_fig} expresses coexistable pairs of two types.
Of course, a parent and its child can coexist.

\section{Proof}

Let $K$ be a knot or link in a bridge position with respect to the standard Morse function $h:S^3\to \Bbb{R}$, and $F$ be a closed incompressible and meridionally incompressible surface.
By Lemma \ref{level}, there exists a level sphere $S$ which intersects $F-K$ essentially and decomposes the pair $(S^3,K)$ into two 3-string trivial tangles $(B_1,T_1)$ and $(B_2,T_2)$.
Put $F_i=F\cap B_i$ for $i=1,2$.

\begin{proof} (of Theorem \ref{I})
Let $F$ be of Type I.
Then, each component of $F_i$ is either of Type $A_0$, $P_0$, $T_0$, $Q_{01}$, $Q_{02}$, $U_{02}$ or $U_{03}$.

Firstly, if $F_1$ has a component of Type $U_{02}$, then $F_1$ consists of one component and $F_2$ consists of one component of Type $A_0$ since Type $U_{02}$ cannot coexist with Type $A_0$.
Hence we have the conclusion I-a.
Similarly, if $F_1$ has a component of Type $U_{03}$, then $F_1$ consists of one component and $F_2$ consists of one component of Type $A_0$.
Hence we have the conclusion I-b.

Secondly, if $F_1$ has a component of Type $Q_{01}$, then $F_1$ consists of one component and $F_2$ consists of two components of Type $A_0$ since Type $Q_{01}$ cannot coexist with Type $A_0$.
However, the boundary of Type $Q_{01}$ does not coincide with the boundary of Type $A_0\cup A_0$, a contradiction.
Similarly, $F_1$ does not have a component of Type $Q_{02}$.

Thirdly, therefore we may assume that each component of $F_1$ and $F_2$ is either $A_0$, $P_0$ or $T_0$ hereafter.
We note that three Types $A_0$, $P_0$ and $T_0$ can coexist.
Let $a$ (resp. $p$, $t$) be  the number of components of $F_1\cup F_2$ whose type is $A_0$ (resp. $P_0$, $T_0$).

If $t=2$, then $p=0$ since $F$ is a genus two closed surface.
We note that two tori of Type $T_0$ are contained in separate tangles, otherwise the boundary does not match.
Thus, $F_1$ has Type $T_0\cup A_0\cup\cdots\cup A_0$ and $F_2$ has Type $T_0\cup A_0\cup\cdots\cup A_0$.
If $a=0$, then we have the conclusion I-a.
If $a=2$, then we have the conclusion I-c.
If $a\ge 4$, then $F$ has more than one component.

If $p=2$, then $t=0$ since $F$ is a genus two closed surface.
We note that two pants of Type $P_0$ are contained in separate tangles, otherwise the boundary does not match.
Thus, $F_1$ has Type $P_0\cup A_0\cup\cdots\cup A_0$ and $F_2$ has Type $P_0\cup A_0\cup\cdots\cup A_0$.
If $a=0$, then we have the conclusion I-b.
If $a=2$, then $F$ has a compressing disk which is formed by gluing two $\partial$-compressing disks for two pants of Type $P_0$ in $B_1$ and $B_2$.
If $a\ge 4$, then $F$ has more than one component.

If $t=1$, then $p=1$ and vice versa.
First suppose that $F_1$ has Type $T_0\cup P_0\cup A_0\cup\cdots\cup A_0$ and $F_2$ has Type $A_0\cup\cdots\cup A_0$.
If $a=2$, namely $F_1$ has Type $T_0\cup P_0$ and $F_2$ has Type $A_0\cup A_0$, then we have the conclusion I-c.
If $a\ge 4$, then the boundary does not match.
Next suppose that $F_1$ has Type $T_0\cup A_0\cup\cdots\cup A_0$ and $F_2$ has Type $P_0\cup A_0\cup\cdots\cup A_0$.
Then, the boundary does not match.
\end{proof}

\begin{proof} (of Theorem \ref{II})
Let $F$ be of Type II.
Then, each component of $F_i$ is either of Type $D_1$, $D_2$, $A_0$, $A_1$, $A_{21}$, $A_{22}$, $A_{23}$, $P_0$, $Q_{01}$, $Q_{02}$, $T_0$, $T_1$, $U_{02}$ or $U_{03}$.

Let $a$ be the number of components of $F_1\cup F_2$ whose Type is $A_0$.

Firstly, if $F_1$ has a component of Type $D_1$, then any component of $F_1$ has Type $D_1$ since Type $D_1$ can coexist with only Type $D_2$.
If $F_1$ is Type $D_1$, then $F_2$ is Type $T_1$ and hence we have the conclusion II-a.
If $F_1$ is Type $D_1\cup D_1$, then $F_2$ is Type $U_{02}$ or $U_{03}$.
In this case, the boundary does not match.

Secondly, if $F_1$ has a component of Type $D_2$, then there are two cases;

\begin{itemize}
	\item [(1)] Each component of $F_1\cup F_2$ has Type $D_2$, $A_0$ or $T_0$.
	\item [(2)] Each component of $F_1\cup F_2$ has Type $D_2$, $A_0$ or $P_0$.
\end{itemize}

In case (1), first suppose that $F_1$ has Type $D_2\cup A_0\cup\cdots\cup A_0$ and $F_2$ has Type $T_0\cup A_0\cup\cdots\cup A_0$.
If $a=0$, then we have the conclusion II-a.
If $a\ge 2$, then the boundary does not match.
Next suppose that $F_1$ has Type $D_2\cup T_0\cup A_0\cup\cdots\cup A_0$ and $F_2$ has Type $A_0\cup\cdots\cup A_0$.
If $a=1$, namely $F_1$ is Type $D_2\cup T_0$ and $F_2$ is Type $A_0$, then we have the conclusion II-c.
If $a\ge 3$, then the boundary does not match.

In case (2), first suppose that $F_1$ has Type $D_2\cup A_0\cup\cdots\cup A_0$ and $F_2$ has Type $P_0\cup A_0\cup\cdots\cup A_0$.
If $a=1$, namely $F_1$ has Type $D_2\cup A_0$ and $F_2$ has Type $P_0$, then we have the conclusion II-c.
If $a\ge 3$, then the boundary does not match.
Next suppose that $F_1$ has Type $D_2\cup P_0\cup A_0\cup\cdots\cup A_0$ and $F_2$ has Type $A_0\cup\cdots\cup A_0$.
In this case, the boundary does not match.

Thirdly, if $F_1$ has a component of Type $A_1$, then any component of $F_1\cup F_2$ has Type $A_0$ or $A_1$.
First suppose that $F_1$ has Type $A_1\cup A_0\cup \cdots \cup A_0$ and $F_2$ has Type $A_1\cup A_0\cup \cdots \cup A_0$.
If $a=0$, then we have the conclusion II-b.
If $a\ge 2$, then $F$ has more than one component.
Next suppose that $F_1$ has Type $A_1\cup A_1\cup A_0\cup \cdots \cup A_0$ and $F_2$ has Type $A_0\cup \cdots \cup A_0$.
In this case, the boundary does not match.

Fourthly, if $F_1$ has a component of Type $A_{21}$, then $F_1$ consists of one component and $F_2$ consists of one component of Type $A_0$ since $A_{21}$ cannot coexist with $A_0$.
Therefore we have the conclusion II-b.

Fifthly, if $F_1$ has a component of Type $A_{22}$, then the boundary of Type $A_{22}$ does not match with the boundary of Type $A_0$.

Sixthly, if $F_1$ has a component of Type $A_{23}$, then $F_1$ consists of one component and $F_2$ consists of one component of Type $A_0$ since $A_{23}$ cannot coexist with $A_0$.
Therefore we have the conclusion II-a.
\end{proof}

\begin{proof} (of Theorem \ref{III})
Let $F$ be of Type III.
Then, each component of $F_i$ is either of Type $D_1$, $D_2$, $D_3$, $A_0$, $A_1$, $A_{21}$, $A_{22}$, $A_{23}$, $P_0$, $Q_{01}$, $Q_{02}$, $P_{11}$ or $P_{12}$.

Firstly, if $F_1$ has a component of Type $D_1$, then there are following cases since Type $D_1$ can coexist with only Type $D_2$.

\begin{enumerate}
	\item $F_1$ is Type $D_1$ and $F_2$ is Type $D_3$.
	\item $F_1$ is Type $D_1\cup D_1$ and $F_2$ is Type $A_{21}$.
	\item $F_1$ is Type $D_1\cup D_1$ and $F_2$ is Type $A_{22}$.
	\item $F_1$ is Type $D_1\cup D_1$ and $F_2$ is Type $A_{23}$.
	\item $F_1$ is Type $D_1\cup D_1\cup D_1$ and $F_2$ is Type $P_{11}$.
	\item $F_1$ is Type $D_1\cup D_1\cup D_1$ and $F_2$ is Type $P_{12}$.
	\item $F_1$ is Type $D_1\cup D_1\cup D_1\cup D_1$ and $F_2$ is Type $Q_{01}$.
	\item $F_1$ is Type $D_1\cup D_1\cup D_1\cup D_1$ and $F_2$ is Type $Q_{02}$.
	\item $F_1$ is Type $D_1\cup D_2$ and $F_2$ is Type $A_{1}$.
	\item $F_1$ is Type $D_1\cup D_1\cup D_2$ and $F_2$ is Type $P_{0}$.
\end{enumerate}

In case 1, we have the conclusion III-a.
In case 3, we have the conclusion III-b.
In case 9, we have the conclusion III-b.
In other cases, the boundary does not match.

Secondly, hereafter we suppose that neither $F_1$ nor $F_2$ has Type $D_1$.
Then, each component of $F_1$ and $F_2$ is Type $D_2$ or $A_0$.
First suppose that $F_1$ is Type $D_2\cup A_0\cup\cdots\cup A_0$ and $F_2$ is Type $D_2\cup A_0\cup\cdots\cup A_0$.
If $a=0$, namyly $F_1$ is Type $D_2$ and $F_2$ is Type $D_2$, then we have the conclusion III-a.
If $a\ge 2$, then $F$ has more than one component.
Next suppose that $F_1$ is Type $D_2\cup D_2\cup A_0\cup\cdots\cup A_0$ and $F_2$ is Type $A_0\cup\cdots\cup A_0$.
If $a=1$, namely $F_1$ is Type $D_2\cup D_2$ and $F_2$ is Type $A_0$, then we have the conclusion III-b.
If $a\ge 3$, then the boundary does not match.
\end{proof}

\section{Future}

Although the number of cases to be distributed is more large, it is possible to characterize closed incompressible surfaces in a given knot complement by the same argument in this paper.
Moreover, since forms and positions of closed incompressible surfaces are restricted when a knot is in bridge position, it is also possible to determine whether a given knot is small or not.

\bigskip

\noindent{\bf Acknowledgement.}
I would like to thank to Prof. Kimihiko Motegi for suggesting to clear the sufficient condition for the surface to be incompressible and meridionally incompressible.

\bibliographystyle{amsplain}

\begin{thebibliography}{10}

\bibitem{HT} A. Hatcher and W. Thurston, {\em Incompressible surfaces in $2$-bridge knot complements}, Invent. Math. {\bf 79} (1985), 225-246.

\bibitem{M} W. Menasco, {\em Closed incompressible surfaces in alternating knot and link complements}, Topology {\bf 23} (1984), 37-44.

\bibitem{O} U. Oertel, {\em Closed incompressible surfaces in complements of star links}, Pac. J. of Math. {\bf 111} (1984), 209-230.

\bibitem{Oz} M. Ozawa, {\em Closed incompressible surfaces in the complements of positive knots}, Comment. Math. Helv. \textbf{77} (2002) 235-243. 


\end{thebibliography}

\end{document}